\documentclass[sn-mathphys-num]{sn-jnl}

\usepackage{graphicx}%
\usepackage{multirow}%
\usepackage{amsmath,amssymb,amsfonts}%
\usepackage{amsthm}%
\usepackage{mathrsfs}%
\usepackage[title]{appendix}%
\usepackage{xcolor}%
\usepackage[dvipsnames]{xcolor}
\usepackage{textcomp}%
\usepackage{manyfoot}%
\usepackage{booktabs}%
\usepackage{algorithm}%
\usepackage{algorithmicx}%
\usepackage{algpseudocode}%
\usepackage{listings}%
\usepackage[all]{xy}
\usepackage{hyperref}

\theoremstyle{thmstyleone}%
\newtheorem{theorem}{Theorem}
\newtheorem{proposition}[theorem]{Proposition}%
\newtheorem{lemma}[theorem]{Lemma}%
\newtheorem{corollary}[theorem]{Corollary}%

\theoremstyle{thmstyletwo}%
\newtheorem{example}{Example}%

\theoremstyle{thmstylethree}%
\newtheorem{definition}{Definition}%

\begin{document}

\title[Some Constructions in the Category of Topological Partial Actions]{Some Constructions in the Category of Topological Partial Actions}


\author*[1]{\fnm{Jesús} \sur{Ávila}}\email{javila@ut.edu.co}
\equalcont{These authors contributed equally to this work.}

\author[2]{\fnm{Fabián} \sur{Molina}}\email{molinafabian@uniovi.es}
\equalcont{These authors contributed equally to this work.}

\affil*[1]{\orgdiv{Departamento de Matem\'aticas y Estad\'istica, Facultad de Ciencias}, \orgname{Universidad del Tolima}, \orgaddress{\street{Cra. 1, C. 42, B. Santa Helena}, \city{Ibagué}, \postcode{CP-730006299}, \state{Tolima}, \country{Colombia}}}

\affil[2]{\orgdiv{Departamento de Matem\'aticas, Facultad de Ciencias}, \orgname{Universidad de Oviedo}, \orgaddress{\street{C. Leopoldo Calvo Sotelo 18}, \city{Oviedo}, \postcode{CP-33007}, \state{Principado de Asturias}, \country{España}}}


\abstract{We present some theoretical aspects of the categories of global and partial group actions of topological groups on topological spaces, denoted by {\bf $G-$ActTop} and {\bf $G-$pActTop}, respectively.  We construct topological partial actions on products and disjoint unions of topological spaces, and $G-$pinvariant subspaces. Furthermore, we provide an adequate definition of congruence ($G-$pcongruence), which enables the construction of quotient objects in the category of topological partial actions. In addition, we associate a special $G-$pcongruence to each pair of strong $G-$pmorphisms. Then, we study products, co-products, equalizers, co-equalizers, pullbacks, and pushouts in the category {\bf $G-$pActTop}. We find that the category of topological partial actions is finitely complete and contains a co-complete subcategory. Finally, we show that the relationship between a topological partial action and its enveloping action is functorial. As a consequence, we obtain that the categorical product in {\bf $G-$ActTop} coincides with the categorical product in {\bf $G-$pActTop}. Moreover, the envelope of a co-product in {\bf $G-$pActTop} is given by the co-product of the respective envelopes.}

\keywords{Topological partial action, category, product, co-product, equalizer, co-equalizer, pullback, pushout, quotient object, adjoint functor.}


\pacs[MSC Classification]{54H15, 18A30, 18A40.}

\maketitle

\section{Introduction}
The study of group actions can be considered to be important for the development of mathematics, since they have been used successfully in almost all its areas \cite{hand1, kerber}. Particularly in algebra, group actions have made it possible to define new classes of rings, such as skew group rings, crossed products, and generalized polynomials, among others.\\

The theory of partial actions has been developed in the last twenty years and it has been shown to be an excellent generalization of global actions. Firstly, Exel in \cite{Exel} defined these actions and demonstrated some of their properties, and later in \cite{Abadie1,Abadie2,Kellendonk} the study of the enveloping action of a partial action was initiated. The importance of this topic is evident in many recent applications in areas such as rings and algebras \cite{Alves-Batista-Dokuchaev-Paques,Bagio-Pinedo, Dokuchaev1,Dokuchaev2,Dokuchaev-Vilca,Ferrero,Jian-Kuo-Szeto}, metric and topological spaces \cite{Abadie2,Choi,Megrelishvili-Schroder,Pinedo-Uzcategui}, homology, cohomology, and homotopy theories \cite{Dokuchaev-Khrypchenko,Dokuchaev-Novikov-Pinedo,Steinberg}, operator theory and dynamical systems \cite{Abadie2,Exel2,Rainone}, among others.\\

In \cite{Ovalle}, the categories of global and partial group actions on sets, denoted by $G-\mathbf{Act}$ and $G-\mathbf{pAct}$ respectively, were defined, and it was proved that the enveloping action of a partial action can be seen as a functor between these categories. Later in \cite{Avila1}, products, co-products, equalizers, pullbacks, coequalizers, and pushouts  in the category $G-\mathbf{pAct}$ were characterized.\\

The purpose of the present paper is to extend several results given in \cite{Avila1,Ovalle} to the category of topological partial actions.
We want to enrich the theory of topological partial actions from the categorical point of view, which will allow us to obtain new results and lines of research. In particular, we study several important concepts in categories, which are contextualized to the category of topological partial actions. In addition, we study, with the functorial approach, the construction of the enveloping action of a topological partial action given in \cite{Abadie2}.
This paper is organised as follows. In Section 2, we present the basic theory of topological global and partial actions and define the categories $G-\mathbf{ActTop}$ and $G-\mathbf{pActTop}$.
Also, we observe that there is a natural inclusion functor $I$ from $G-\mathbf{ActTop}$ to $G-\mathbf{pActTop}$. In Section 3, we study in detail several concepts in the category $G-\mathbf{pActTop}$.
We characterize their finite products, co-products, equalizers, pullbacks, and construct quotient objects. In particular, we show that the category  $G-\mathbf{pActTop}$ is finitely complete. Moreover, we define the category $G-\mathbf{spActTop}$, which is a subcategory of $G-\mathbf{pActTop}$ and study its co-equalizers and pushouts. In particular, we prove that the category $G-\mathbf{spActTop}$ is co-complete. Finally, in Section 4 the construction of the enveloping action given in \cite{Abadie2} is presented in a general way as the minimal and unique globalization of a topological partial action. Moreover, we prove that this construction of the enveloping action can be seen as a covariant functor $E$ from $G-\mathbf{pActTop}$ to $G-\mathbf{ActTop}$, which is called the enveloping functor, and show that it is a left adjoint of the inclusion functor $I$. Then, we apply that adjunction to prove that the categorical product in $G-\mathbf{ActTop}$ coincides with the respective categorical product in $G-\mathbf{pActTop}$, and that the envelope of a disjoint union coincides with the disjoint union of the envelopes.

\section{Categories of Topological Actions}

In this section we integrate the theory of topological global actions with the theory of categories. For this, it is necessary to define the following notions.\\

\begin{definition}
Let $ G $ be a topological group with identity element $ e $ and $ X $  a topological space. A topological global  action of $ G $ on $ X $ is a function $ \varphi: G \times X \to X $ that satisfies:
\begin{enumerate}
  \item $ \varphi (e, x) = x $ for each $ x \in X $.
  \item $ \varphi (gh, x) = \varphi (g, \varphi (h, x)) $, for each pair $ g, h \in G $ and each $ x \in X $.
  \item $ \varphi $ is continuous with the product topology on $ G \times X $.
\end{enumerate}
In this case, we say that $ G $ acts globally and topologically on $ X $.\\
\end{definition}

Note that if $ \varphi $ is a topological global  action of $ G $ on $ X $, then for each $ g \in G $, the function $ \varphi_g: X \to X $ defined by $ \varphi_g (x) = \varphi (g, x) $ is a homeomorphism.
In fact, $ \varphi_g $ is bijective for each $ g \in G$.
Now, let $ x \in X $, $ g \in G $ and $ U \subseteq X $ be an open neighborhood of $ \varphi_g (x) $. Then $ \varphi (g, x) \in U$.
Since $ \varphi $ is continuous, there are open sets $ H \subseteq G $ and $ V \subseteq X $ such that $ (g, x) \in H \times V $ and $ \varphi (H \times V) \subseteq U $.
Thus, $ V $ is an open neighborhood of $ x $ such that $ \varphi_g (V) = \varphi (\{g \} \times V) \subseteq \varphi (H \times V) \subseteq U $. Therefore, $ \varphi_g $ is continuous.
On the other hand,  $ \varphi_g^{-1} = \varphi_ {g^{- 1}} $ is continuous.\\

\begin{definition}[\cite{Abadie2}]\label{defpartop}
Let $ G $ be a topological group with identity element $ e $ and $ X $  a topological space. It is said that $ \alpha $ is a topological  partial action of $ G $ on $ X $, if there exists a class of pairs $ \{(X_g, \alpha_g) \}_{g \in G}$, where $ X_g $ is an open subset of $ X $ and $ \alpha_g: X_ {g^{- 1}} \to X_g $ is a homeomorphism, such that:
\begin{enumerate}
\item $X_e = X$ y $\alpha_e=id_X$.
\item $\alpha^{-1}_h(X_h\cap X_{g^{-1}})\subseteq X_{(gh)^{-1}}$ for each pair $ g, h \in G $
\item $(\alpha_g\circ \alpha_h)(x) = \alpha_{gh}(x)$  for each  $x\in \alpha^{-1}_h (X_h\cap X_{g^{-1}})$.
\item The set $\Gamma_\alpha=\{(g,x)\in G\times X: x\in X_{g^{-1}}\}$ is open in $G \times X$.
\item The function $\alpha:\Gamma_\alpha\to X$  given by $\alpha(g,x)=\alpha_g(x)$,  is continuous.
  \end{enumerate}
In this case, we also say that $ G $ acts partially and topologically on $ X $.\\
\end{definition}

The following proposition states that the restriction of a topological global action to an open set of a topological space is a topological partial  action.
This result is also important because it gives rise to the notion of the globalization of a topological partial action, which will be studied in the last section.\\

\begin{proposition}[\cite{Abadie2}]
Let $ \varphi $ be  a topological global  action of $ G $ on $ X $. Let $ S $ be an open subset of $ X $, $ \{(S_g, \alpha_g) \} _ {g \in G} $ the class given by $ S_g = S \cap \varphi_g (S) $, and $ \alpha_g: S_ {g^{- 1}} \to S_g $  defined as $ \alpha_g (x) = \varphi_g (x) $ for each $ x \in S_ {g^{- 1}} $. Then the function $ \alpha: \Gamma _ {\alpha} \to S $ given by $ \alpha (g, x) = \alpha_{g} (x) $ where $ \Gamma_ {\alpha} = \{(g, x) \in G \times S: x \in S_{g^{-1}}\}$, is a topological partial  action of $ G $ on $ S $. This restriction of $ \varphi $ is called the induced topological partial action of $ G $ on $ S $.\\
\end{proposition}
From the previous proposition, it is clear that every topological global action is a  topological partial action.
In fact, if $ \varphi $ is a topological global  action of $ G $ on $ X $, we can take $S=X$ and then $ \{(X, \varphi_g) \}_{g \in G} $ is a class of pairs that satisfy conditions 1--5 of Definition \ref{defpartop}.\\

On the other hand, if $ \alpha $ is a topological partial action of $ G $ on $ X $, there exists a class of pairs $ \{(X_g, \alpha_g) \}_{g \in G} $ such that:
\begin{enumerate}
\item For each $ x \in X $, we have that $ (e, x) \in \Gamma_ \alpha $ since $X = X_e $. Hence, $ \alpha (e, x) = \alpha_e (x) = x $ for each $ x \in X $.
\item If $ x \in \alpha^{- 1}_h (X_{h} \cap X_{g ^ {- 1}})$, then $ \alpha_h (x) \in X_{h} \cap X_{g^{- 1}} $. Thus:
\begin{enumerate}
  \item $ x \in X_{h^{- 1}}$ and thus $(h, x)\in \Gamma_\alpha $.
  \item $ \alpha_h (x) \in X_{g^{- 1}}$ and so  $ (g, \alpha_h(x)) \in \Gamma_ \alpha $.
\end{enumerate}
Moreover, $ x \in X_{(gh)^{- 1}}$ and thus $ (gh, x) \in$ $\Gamma_\alpha $. Therefore,
\begin{align*}
\alpha (gh, x) &= \alpha_ {gh} (x) = \alpha_g (\alpha_h (x))= \alpha (g, \alpha (h, x)).
 \end{align*}
\item The function $ \alpha: \Gamma_ \alpha \to X $ given by $ \alpha (g, x) = \alpha_g (x) $ is continuous.
  \end{enumerate}
Thus, $\alpha$ is a topological global action of $G$ on $X$ if, and only if,
$\Gamma_{\alpha} = G \times X$. In what follows, we present some examples of
topological partial actions that do not satisfy this condition, that is,
they are not topological global actions.\\

\begin{example}[\cite{Abadie2}]\label{ejempartop2}
Let $ X $ be a topological space and let $ h: X \to X $ be a homeomorphism. If $ \{(X_t, \alpha_t) \}_{t \in \mathbb{R}} $ is given by
$ X_t = X $ and $ \alpha_t = h^t $ when $ t \in \mathbb{Z} $, and $ X_t = \emptyset $ and $ \alpha_t = \emptyset $ when $ t \not \in \mathbb{Z} $. Then,  the function $ \alpha: \Gamma_ \alpha \to X $ given by $ \alpha (t, x) = h^t (x) $, is a topological partial  action of  $\mathbb{R}_d$ (with the discrete topology) on  $ X $ where $ \Gamma_ \alpha = \mathbb{Z} \times X $. 
\\
\end{example}

\begin{example}[\cite{Abadie2}]\label{ejempartop2b}
Suppose that $G=\{1,-1\}$ endowed with the discrete topology, $X=[0,1]$ endowed with the usual subspace topology, and $\alpha$ is the topological partial action of $G$ on $X$ given by
\(
X_1 = X\),  \(X_{-1} = V\),  \( \alpha_1 = i_X\),  \(\alpha_{-1} = i_V,
\)
where $V=(a,1]$ for some fixed $0<a<1$. We have $X_{1^{-1}}=X_1=X=[0,1]$ and $X_{(-1)^{-1}}=X_{-1}=V=(a,1]$, hence $(1,x)\in\Gamma_\alpha$ for all $x\in[0,1]$ and  $(-1,x)\in\Gamma_\alpha$ for all $x\in(a,1]$. Therefore,
\(
\Gamma_\alpha
=
\big(\{1\}\times[0,1]\big)
\cup
\big(\{-1\}\times(a,1]\big)
\). \\
\end{example}

\begin{example}[\cite{Choi}]
Let $X=\mathbb{R}_d$ and let $G=GL(2,\mathbb{R})$ endowed with the usual topology.
For each
$g=\begin{pmatrix} a & b \\ c & d \end{pmatrix}\in G$, define
\(
X_g=\{x\in\mathbb{R}:cx+d\neq 0\}
\)
and let $\alpha_g:X_{g^{-1}}\to X_g$ be the function given by
\[
\alpha_g(x)=\frac{ax+b}{cx+d}.
\]
Then the family $\{(X_g,\alpha_g)\}_{g\in G}$ defines a topological partial
action of $G$ on $X$. In this case,
\(
\Gamma_\alpha=\{(g,x)\in GL(2,\mathbb{R})\times\mathbb{R}:c'x+d'\neq 0\},
\)
where
$g^{-1}=\begin{pmatrix} a' & b' \\ c' & d' \end{pmatrix}$.\\
\end{example}

\begin{example}[\cite{Avila-Molina}]

Let $X=\mathbb{R}$ endowed with the usual topology and consider the family
$\{(X_t,\alpha_t)\}_{t\in\mathbb{R}}$, where
$X_t=\{x\in\mathbb{R}:t<e^x\}$ and $\alpha_t:X_{-t}\to X_t$ is defined by
$\alpha_t(x)=\ln(t+e^x)$. Then the map $\alpha:\Gamma_\alpha\to X$ given by
$\alpha(t,x)=\alpha_t(x)$ defines a topological partial action of $\mathbb{R}$
on $X$, where $\Gamma_\alpha=\{(t,x)\in\mathbb{R}^2:-t<e^x\}$. 
The previous example is a particular case of the flow of a differentiable vector field viewed as a
topological partial action of $\mathbb{R}$ on a differentiable manifold.\\
\end{example}

\begin{example}[\cite{Abadie2}]\label{ejempartop3}

Let $X$ be a $C^2$ differentiable manifold and let
$v : X \to TX$ be a $C^1$ vector field on $X$. For each $x \in X$, denote by
$\gamma_x$ the integral curve of $v$ satisfying $\gamma_x(0)=x$ and
$\gamma_x'(t)=v(\gamma_x(t))$, whose domain is a maximal open interval
$I_x \subseteq \mathbb{R}$ containing $0$. For each $t \in \mathbb{R}$, define
\(
X_{-t}=\{x\in X : t\in I_x\}
\quad \text{and} \quad
\alpha_t : X_{-t} \to X_t\)
by $\alpha_t(x)=\gamma_x(t)$. Then
\(
\Gamma_\alpha=\{(t,x)\in \mathbb{R}\times X : -t\in I_x\},
\)
and consequently the map $\alpha : \Gamma_\alpha \to X$ given by
$\alpha(t,x)=\alpha_t(x)$ defines a topological partial action of
$\mathbb{R}$ on $X$. The details of the proof can be found in \cite{Avila-Molina}.\\
\end{example}

With some conditions on the topological group and on the topological space, it is possible to determine when a topological partial action is a topological global action.\\

\begin{theorem}[\cite{Abadie2}]
If $ \alpha $ is a topological  partial action of $ G $ on a compact topological space $ X $, there exists an open subgroup $ H \subseteq G $ such that $ H \times X \subseteq \Gamma_\alpha $ and the restriction of $ \alpha $ to $ H \times X $ is a topological global  action of $ H $ on $ X $. Moreover, if $ G $ is connected then $ \alpha $ is a topological  global action of $ G $ on $ X $.\\
\end{theorem}

As one of the objectives of this paper is to integrate the theory of topological global actions  with the theory of categories,  it is necessary to define the following notion.\\

\begin{definition}
Let $ \varphi $ and $ \psi $ be topological global  actions of $ G $ on $ X $ and $ Y $ respectively. We say that a continuous function $ f: X \to Y $ is a $ G-$morphism of $ X $ in $ Y $ (or simply $ G- $morphism), if  $ f (\varphi (g, x)) = \psi (g, f (x)) $ for each pair $ x \in X $ and $ g \in G $.\\
\end{definition}

It is easy to prove that given  a  topological global action  of $ G $ on $ X $,  the identity function $ i_{X} $ is a $G-$morphism  and that the composition of $G-$morphisms is also a $G-$morphism.\\

Thus, given a topological group $G$, it is possible to define the category of topological global actions of $G$, denoted by $G-\mathbf{ActTop}$, as follows:
\begin{enumerate}
\item $\text{Ob}(G-\mathbf{ActTop}) $ is the class of all ordered pairs $(X,\varphi)$ where $X$ is a topological space and $\varphi$ is a  topological global action of $G$ on  $X$.  So, every object $(X,\varphi)$ of   $G-\mathbf{ActTop}$ will be denoted by $X$, by assuming that there is an implicit topological global action of $G$ on $X$.
\item For each pair $ X, Y \in \text{Ob} (G-\mathbf{ActTop}) $, $ \text{Mor}(X, Y) $ is the set of all\linebreak $G-$morphisms of $ X $ to $ Y $.
\item The internal composition law $ \circ $ is the usual composition of functions. This satisfies:
  \begin{enumerate}
\item If $ f_1 \in \text{Mor} (X, Y) $ and $ f_2 \in \text{Mor} (Y, Z) $, then $ f_2 \circ f_1 \in \text {Mor} (X , Z) $.
\item Since $ \circ $ is associative, the associativity of the $ G- $ morphisms is guaranteed.
\item If $ X \in \text{Ob} (G-\mathbf{ActTop}) $, then $ i_{X} \in \text{Mor} (X, X) $.
  \end{enumerate}
  \end{enumerate}

\begin{definition}[\cite{Abadie2}]\label{defmor}
Let  $ \alpha $ and $ \beta $ be topological partial actions of $ G $ on $ X $ and $ Y $ respectively. A continuous function $ f: X \to Y $ is a $ G-$pmorphism from $ X $ to $ Y $(or simply $G-$pmorphism), if $ f (X_g) \subseteq Y_g $ and $ f (\alpha_g (x)) = \beta_g (f (x)) $, for each pair $ g \in G $ and $ x \in X_{g^{- 1}}$.\\
\end{definition}

From this definition, we can also verify that given a topological partial action  of $ G $ on $ X $, the identity function $ i_ {X} $ is a $ G-$pmorphism and that the composition of $G-$pmorphisms is also a $G-$pmorphism.\\

Hence,  given a topological group $ G $, it is possible to define the category of topological partial actions of $G$, denoted by $G-\mathbf{pActTop}$, as follows:
\begin{enumerate}
\item $\text{Ob} (G-\mathbf{pActTop}) $ is the class of all ordered pairs $(X,\alpha)$ where $X$ is a topological space and $\alpha$ is a  topological partial action of $G$ on  $X$. Thus, every object $(X,\alpha)$  of $G-\mathbf{pActTop}$ will be denoted by $X$, by assuming that there is an implicit topological partial action of $G$ on $X$.
\item If $ X, Y \in \text{Ob} (G-\mathbf{pActTop}) $, then $ \text{Mor}(X, Y) $ is the set of all $ G-$pmorphisms of $ X $ to $ Y $.
\item The internal composition law is the usual composition $ \circ $ of functions. This satisfies:
  \begin{enumerate}
\item If $ f_1 \in \text{Mor} (X, Y) $ and $ f_2 \in \text{Mor} (Y, Z) $, then $ f_2 \circ f_1 \in \text{Mor} (X , Z) $.
\item Since $ \circ $ is associative, the associativity of the $G-$pmorphisms is guaranteed.
\item If $ X \in \text{Ob} (G-\mathbf{pActTop})$, we have $ i_ {X} \in \text{Mor} (X, X) $.
  \end{enumerate}
  \end{enumerate}

Note that since every topological global action is a topological partial action, it is natural to define $I:(G-\mathbf{ActTop}) \to(G-\mathbf{pActTop}) $ in the following way and show that it is a covariant functor:
\begin{enumerate}
\item If $ X \in \text{Ob} (G-\mathbf{ActTop}) $, then $ I (X) = X \in \text{Ob}(G-\mathbf{pActTop}) $.
\item If $ X, Y \in \text{Ob} (G-\mathbf{ActTop}) $ and $ f \in \text{Mor} (X, Y) $, then $ I (f ): I (X) \to I (Y) $ given by $ I (f) = f $ is an element of $ \text{Mor} (I (X), I (Y)) $. In fact, if $ \alpha $ and $ \beta $ are topological global actions of $ G $ on $ X $ and $ Y $ respectively, then $ \alpha $ and $ \beta $ are topological partial actions of $ G $ on $ X $ and $ Y $. Thus, $$ f (X_g) = f (X) \subseteq Y = Y_g,$$
for all $g\in G$. Also, for each $ x \in X_{g^{-1}}$, we have that $ f (x) \in Y_{g^{-1}} $ and so
      \begin{align*}
      f (\alpha_g (x)) &= f (\alpha (g, x)) = \beta (g, f (x)) = \beta_g (f (x)).
      \end{align*}
      That is, $ I (f) = f $ is a $G-$pmorphism.
  \item If $ X, Y, Z \in \text{Ob} (G-\mathbf{ActTop}) $, $ f_1 \in \text {Mor} (X, Y) $, $ f_2 \in \text {Mor} (Y, Z) $ and $ x \in X $,
   we have
  $$ I (f_2 \circ f_1) (x) = (f_2 \circ f_1) (x) = (I (f_2) \circ I (f_1)) (x).$$
  \item For each $ x \in X $ we have $ I (i_{X}) (x) = i_{X} (x) = i_{I (X)} (x) $.
\end{enumerate}
This functor is the inclusion from $G-\mathbf{ActTop}$ to $G-\mathbf{pActTop}$. This means that the category  $G-\mathbf{ActTop}$ can be seen as a full subcategory of $G-\mathbf{pActTop}$. With this idea in mind, in the next section we will characterize some categorical notions and demonstrate some interesting properties of the category $G-\mathbf{pActTop}$, which are of course valid in $G-\mathbf{ActTop}$.

\section{Categorical Aspects of $G-\mathbf{pActTop}$}

In this section we study the category $G-\mathbf{pActTop}$. We prove that this category has finite products, co-products, equalizers, pullbacks, and quotient objects. In particular, we obtain  that $G-\mathbf{pActTop}$ is finitely complete. Moreover, we define the category $G-\mathbf{spActTop}$, which is a subcategory of $G-\mathbf{pActTop}$ and study its co-equalizers and pushouts. In particular, we prove that the category $G-\mathbf{spActTop}$ is co-complete. The categorical concepts used in this paper can be found in \cite{Adamek,MacLane}.

\subsection{Study of some limits}
Given  a collection $\{X_i\}_{i \in I} $ of topological spaces, if $\alpha_i$ is a topological partial action of $G$ on $X_i$ for all $i\in I$, the most natural way to define a topological partial action of $G$ on $ X = \prod_{i \in I} X_i $, with the product topology, is taking  $X_g$ as  $\prod_{i \in I} X_{i,g}$ and $ \alpha: \Gamma_{\alpha} \to X $ as  $ \alpha(g,x) = (\alpha_i(g, x))_{i\in I}$, for all   $(g,x)\in \Gamma_{\alpha}$.  This is used in the category $G-\mathbf{pAct}$ (see \cite{Avila1}) and is essential to demonstrate that this category is complete. However, when considering topological partial actions, the set $\prod_{i \in I} X_{i,g}$ in that topology, is not always open in $X$ nor is $ \alpha: \Gamma_{\alpha} \to X $ always continuous. Therefore, it is necessary to consider the box topology on $ X = \prod_{i \in I} X_i $. The definitions of these topologies can be found in \cite{munkres}. \\
\begin{proposition}\label{PartialActionPro}
Let $ \{X_i \}_{i \in I} $ be a collection of topological spaces and $ X = \prod_{i \in I} X_i $. If $ \alpha_i $ is a topological partial  action of $ G $ on $ X_i$ for every $i\in I$, then the function $ \alpha: \Gamma_\alpha \to X $ given by $ \alpha(g, x) = (\alpha_i (g , x_i))_{i \in I} $, for all 
$ (g, x) \in \Gamma_{\alpha}$, where $\Gamma_{\alpha} = \{(h, y) \in G \times X: y \in X_{h^{-1}}\}$, is a topological  partial action of $ G $ on $ X$. 
\end{proposition}
\begin{proof}
By assumption, for every $i \in I $, there exists a class $ \{(X_ {i, g}, \alpha_{i, g}) \}_{g \in G} $ such that for each $ g \in G $, $ X_{i, g} \subseteq X_i $ is an open set, $ \alpha_ {i, g}: X_{i, g^{-1}} \to X_{i, g} $ is a homeomorphism, the set $ \Gamma_{\alpha_i} = \{ (g, x_i) \in G \times X_i: x_i \in X_{i,g^{- 1}} \}$ is open in $ G \times X_i $ and the function $ \alpha_i: \Gamma_{\alpha_i} \to X $ given by $ \alpha_i (g, x_i) = \alpha_{i, g} (x_i) $ is continuous and is also a partial action  of $ G $ on $ X_i $.
Now, for each $ g \in G $ let us consider $ X_ {g} = \prod_{i \in I} X_ {i, g} $ and the function $ \alpha_ {g}: X_{g^{-1}} \to X_{g} $ given by $ \alpha_g(x) = (\alpha_{i, g} (x_i))_{i \in I} $ for each $ x = (x_i)_{i \in I} \in X $.
Since $ X_ {i, g} $ is an open set for every  $ i \in I $, then so is $ X_g $ for each $ g \in G $ (with the box topology).
Also, since $ \alpha_{i, g} $ is a homeomorphism for every $ i \in I $, then $ \alpha_g $ is a homeomorphism for each $ g \in G $.\\

Again  items 1, 2 and 3 of Definition \ref{defpartop} are true (\cite{Avila1}). We need only to prove that $\Gamma_{\alpha}$ is open and that the corresponding function $ \alpha: \Gamma_{\alpha} \to X $ is well defined and is continuous.\\

If $ (g, x) \in \Gamma_{\alpha} = \{(h, y) \in G \times X: y \in X_{h^{- 1}}\}$, then $ x \in X_{g^{- 1}} $. Therefore, $ x_i \in X_ {i, g^{- 1}} $  and so $ (g, x_i)  \in \Gamma_{\alpha_{i}} $ for every  $ i \in I $. If $ i_0 \in I $, then $ (g, x_ {i_0}) \in \Gamma_{\alpha_{i_0}} $ and therefore $ (g, x) \in \Gamma_{\alpha_{i_0 }} \times \prod_{i \neq i_0} X_{i, g} = V $. Since $ X_{i, g^{- 1}} $ is open in $ X_i $ and $\Gamma_{\alpha_{i_0}}$ is open in $ G \times X_{i_0} $, it follows that $ V $ is open in $ G \times X $ and in addition $ V \subseteq \Gamma_\alpha $. Hence, $ \Gamma_{\alpha}$  is open in $ G \times X $.\\

Let us notice that $ \alpha (g, x)  = (\alpha_i (g, x_i))_{i\in I} = (\alpha_ {i, g}(x_i))_{i \in I}= \alpha_{g} (x)$, for all $x\in X_{g^{-1}}$, is well defined. Suppose that $ (g, x) \in \Gamma_ \alpha $ and $V\subseteq X$ is an open neighborhood of $\alpha(g,x)$. Then  $ x \in X_{g^{- 1}} $ and $ V= \prod_{i \in I} V_{i}$ where $V_i$ is open in $X_i$ for all $i\in I$. So, $ x_i\in X_{i, g^{- 1 }} $ and $ x_i \in V_i $ for every  $ i \in I $, and since  the function $ \alpha_i: \Gamma _ {\alpha_i} \to X_i $  is continuous, there exist an open neighborhood $U_i\subseteq X_i$ of $x_i$ such that $\alpha_i(U_i)\subseteq V_i$. Let us notice that, if $i_0\in I$, there exist open sets $M\subseteq G$ and $N\subseteq X_{i_0}$ such that $M\times N$ is an open neighborhood of $(g,x_{i_0})$ contained in $U_{i_0}$.  That occurs because $\Gamma_{i_0}$ is open in $G\times X_{i_0}$.  Then, $\alpha_{i_0}(g,x_{i_0})\in \alpha_{i_0}(M\times N)\subseteq \alpha_{i_0}(U_{i_0})$. Consequently, 
$U=M\times N \times \prod_{i \neq i_0} U_{i, g}$ is open in $\Gamma_\alpha$, $(g,x)\in U$ and $\alpha (U)=\alpha_{i_0}(M\times N) \times \prod_{i \neq i_0} \alpha_i(U_{i})\subseteq \prod_{i \in I} V_{i}=V.$
Thus, $ \alpha: \Gamma _ {\alpha} \to X $ is continuous.
\end{proof}

Unfortunately, the topological partial action of  Proposition \ref{PartialActionPro} cannot be characterized as the product of the collection $ \{X_i \}_{i \in I} $ (with the box topology) in the category $G-\mathbf{pActTop}$. 
To see this fact,  let us keep the notation of such a proposition and let $ \pi_i: X \to X_i $ be the projection of $ X = \prod_{i \in I}X_i$ on $ X_i $. Then $\pi_i$ is continuous, $\pi_i(X_g)=X_{i,g}$ and $ \pi_i (\alpha_g (x)) = \alpha_{i,g}(x_i) = \alpha_{i,g}(\pi_i (x)) $ for each pair $ g \in G $ and $ x = (x_i)_{i \in I} \in X_{g^{-1}}$. Therefore, $ \pi_i $ is a $G-$pmorphism.\\

Assuming that $ \beta $ is a partial topological action of $G$ on $Y$ and $ f_i: Y \to X_i $ is a $G-$pmorphism for every $ i \in I$, we have that $ f_i $ is continuous, $f_i(Y_g)\subseteq X_{i,g}$ and $ f_i (\beta_g(y)) = \alpha_{i,g}(f_i (y)) $ for each pair $ g \in G $ and $y \in Y_{g^{-1}} $. Now, if $ y \in Y $, then $ f_i (y) \in X_i $ and so the function $ f: Y \to X $ given by $ f (y) = (f_i (y))_{i \in I} $ is well defined and satisfies that
$$f(Y_g)=\prod_{i\in I}f_i(Y_g)\subseteq \prod_{i\in I}X_{i,g}=X_g$$
and
\begin{align*}
f (\beta_g(y)) = (f_i (\beta_g(y)))_{i \in I} = (\alpha_i,g(f_i (y)))_{i \in I}
= \alpha_g(f (y)),
\end{align*}
for each pair $ g \in G $ and $ y \in Y_{g^{-1}} $. But $ f $ is not a $G-$pmorphism since it  is not necessarily continuous with the box topology on $ \prod_{i \in I}X_i$.\\

However, if $I$ is a finite set and we assume all of the above, then the box topology  coincides with the product topology on  $ \prod_{i \in I}X_i$ and the function $ f: Y \to X $ is continuous. Indeed, if $y\in Y$ and $U\subseteq X$ is an open neighborhood  of $f(y)$, then $U=\prod_{i \in I}U_i$, where $U_i$ is open in $X_i$  and $f_i(y)\in U_i$ for all $i\in I$. Since 
$f_{i}:Y \to X_i $ is continuous, then there exists $V_{i}$ open in $Y$ such that $y\in V_{i}$ and $f_{i}(V_{i})\subseteq U_i$. Consequently, $V=\cap_{i \in I}V_i$ is an open neighborhood of $y$ (because $I$ is finite) and $$f(V)=  \prod_{i \in I} f_i(V) \subseteq \prod_{i \in I} f_i(V_i) \subseteq  \prod_{i \in I}  U_i=U.$$
In that case, $ f $ is the unique function such that $ \pi_i \circ f = f_i $ for every $ i \in I $. Indeed,  $ f_i (y) = \pi_i (f (y)) = (\pi_i \circ f) (y) $ for each $ y \in Y $, that is, the diagram below is commutative:
$$
\xymatrix@=1.5cm{
 & \prod_{i \in I}X_i\ar[d]^{\pi_i}\\
Y \ar@{-->}[ur]^f \ar[r]_{f_i}&X_i
}
$$
Also, if $ h: Y \to X $ is a $G-$pmorphism such that $ \pi_i \circ h = f_i $ for every $ i \in I $, then $ f (y) = (f_i (y))_{i \in I} = (\pi_i(h (y)))_{ i \in I} = h (y) $ for all $ y \in Y $. That is, $ h = f $. Thus, we have obtained the following result.\\

\begin{theorem}\label{charprod}
Let $ \{X_i \}_{i \in I} $ be a finite collection of topological spaces. Then the topological  partial action given by Proposition \ref{PartialActionPro} on 
  $ X = \prod_{i \in I} X_i $ characterizes the product of the collection $ \{X_i \}_{i \in I} $ (with the box or product topology) in the category $G-\mathbf{pActTop}$.\\
\end{theorem}

Although it was not possible to characterize the product object in the category $G-\mathbf{pActTop}$, it is possible to do so in the category $G-\mathbf{ActTop}$, as shown below.\\

\begin{theorem}\label{proglo}
Let $ \{X_i \}_{i \in I} $ be a collection of topological spaces. If $ \varphi_i $ is a global topological action of $ G $ on $ X_i $ for every $i\in I$, then the function $ \varphi: G \times X \to X $ given by $ \varphi (g, x) = (\varphi_i (g , x_i))_{i \in I},$ is a global topological action of $ G $ on $ X = \prod_{i \in I} X_i $. Furthermore, this topological global action characterizes the product of the collection $ \{X_i \}_{i \in I} $ (with the product topology) in the category $ G- \mathbf{ActTop}$.
\end{theorem}
\begin{proof}
If $ x = (x_i)_{i \in I} \in X $, then $ \varphi (e, x) = (\varphi_i(e, x_i ))_{i \in I} = (x_i)_{ i \in I} = x.$ Note that if $ g, h \in G $ and $ x = (x_i)_{ i \in I} \in X $, then $ \varphi (gh, x) = (\varphi_i (gh, x_i))_{i \in I} = (\varphi_i (g, \varphi_i (h, x_i)))_{i \in I} = \varphi (g, \varphi (h, x)) $.  Now, $\varphi_i$ is continuous for every $ i \in I $, then $ \varphi $ is continuous and therefore, is a global topological action of $ G $ on $ X $. \\

Moreover, if $ \pi_i: X \to X_i $ is the projection of $ X $ on $ X_i $, then $\pi_i$ is continuous and $ \pi_i (\varphi (g, x)) = \varphi_i(g, x_i) = \varphi_i(g, \pi_i (x)) $, for each pair $ g \in G $ and $ x = (x_i)_{i \in I} \in X $. So $ \pi_i $ is a $ G-$morphism. \\

If  $ \psi$ is a topological global action of $G$ on $Y$ and $ f_i: Y \to X_i $ is a $ G-$morphism for every $ i \in I$, then $ f_i $ is continuous and $ f_i (\psi (g, y)) = \varphi_i (g, f_i (y)) $, for each pair $ g \in G $ and $y \in Y $. In addition, if $ y \in Y $, then $ f_i (y) \in X_i $ and so the function $ f: Y \to X $ given by $ f (y) = (f_i (y))_{i \in I} $, is well defined and is continuous. Furthemore $$ f (\psi (g, y)) = (f_i (\psi (g, y)))_{i \in I} = (\varphi_i (g, f_i (y)))_{i \in I} = \varphi (g, f (y)), $$ for each pair $ g \in G $ and $ y \in Y $. This implies that $ f $ is a $ G-$morphism. \\

Finally, $ f $ is the unique function such that $ \pi_i \circ f = f_i $ for every $ i \in I $. Indeed,  $ f_i (y) = \pi_i (f (y)) = (\pi_i \circ f) (y) $ for each $ y \in Y$. Also, if $ h: Y \to X $ is a $ G- $morphism such that $ \pi_i \circ h = f_i $ for every $ i \in I $, then $ f (y) = (f_i (y))_{i \in I} = (\pi_i(h (y)))_{ i \in I} = h (y) $ for each $ y \in Y $. Hence $ h = f $.
$$
\xymatrix@=1.5cm{
 & \prod_{i \in I}X_i\ar[d]^{\pi_i}\\
Y \ar@{-->}[ur]^f \ar[r]_{f_i}&X_i
}
$$
\end{proof}


The result given in Theorem \ref{charprod} is important to characterize pullbacks in $G-\mathbf{pActTop}$. The following notion also plays a fundamental role.\\

\begin{definition}
Let $ \alpha $ be a topological partial  action of $ G $ on $ X $. We say that $ A \subseteq X $ is $ G-$pinvariant (with respect to $ \alpha $), if $ \alpha_g (A \cap X_ {g^{-1}}) = A \cap X_g $ for each $ g \in G $.\\
\end{definition}
We  claim that if $ \alpha_g (A \cap X_{g^{- 1}}) \subseteq A \cap X_g $ for each $ g \in G $, then $ A $ is $G-$pinvariant. To prove this fact, just take $g=h^{-1}$ for any $h\in G$ and substitute it in the previous expression.
Furthermore, if $ A $ is $ G-$pinvariant, the restriction of $ \alpha $ to $ A $ is a topological partial action of $ G $ on $ A $.
Indeed, let  $ \{(S_g, \delta_g) \} _ {g \in G} $ be  the class   given by $ S_g = A \cap X_g  $ and $ \delta_g: S_ {g^{- 1}} \to S_g $  defined as $ \delta_g (x) = \alpha_g (x) $ for each $ x \in S_ {g^{- 1}} $.
Now, for each $ g \in G $, $ X_g  $ is an open set. Hence, $ S_g $ is an open set in $A$.
Moreover,  $ \alpha_g: X_{g^{- 1}} \to X_g $ is a homeomorphism  and $ \alpha_g (A \cap X_ {g^{-1}}) = A \cap X_g $ for each $ g \in G $.
Therefore, $ \delta_g: S_{g^{- 1}} \to S_g $  is also a homeomorphism.\\

Since items 1, 2 and 3 of Definition \ref{defpartop} are evident (see \cite{Avila1}), we must to prove that the set $\Gamma_{\delta}$ is open and that the corresponding function $ \delta: \Gamma_{\delta} \to A $ is continuous.\\

If  $ (g, x) \in \Gamma_\delta=\{(h,y)\in G\times A: y \in S_{h^{- 1}}\} $, then $ x \in  A\cap X_{g^{-1}}$. Hence, $ (g, x) \in \Gamma_\alpha $ and thus $ \Gamma_\delta\subseteq \Gamma_\alpha$. Since
    $\Gamma_\alpha$ is open, there are open sets $ U \subseteq G $ and $ V \subseteq X $  such that $ (g, x) \in U \times V $ and $ U \times V \subseteq \Gamma_\alpha $. Consequently, $ U $ and $ A\cap V $ are open sets such that $ (g, x) \in U \times (A\cap V) $ and $ U \times (A\cap V) \subseteq \Gamma_\delta $. In fact, for each $ (h,y) \in U \times (A\cap V)  $ we have that $ (h,y) \in U \times V\subseteq  \Gamma_\alpha $ and thus $ y \in X_{h^{-1}}$. Hence, $ y \in A\cap X_{h^{-1}}=S_{h^{-1}}$  and therefore $ (h, y) \in \Gamma_\delta $.\\

Let $ \delta: \Gamma_\delta \to A $ be the function given by $ \delta(g, x) = \delta_g (x) $. If $ (g, x) \in \Gamma_\delta $, then  $ x\in A\cap X_{g^{-1}}$ and $ \delta(g, x) = \delta_g (x)=\alpha_g (x) = \alpha(g, x)$. Consequently $\delta$ is the restriction of $\alpha$ to $\Gamma_\delta$. Since $\alpha$ is continuous and $\Gamma_\delta $ is open in $G\times A$, then  $\delta$ is continuous.\\

The following results extend Theorem 8 and Corollary 9 of \cite{Avila1} to the topological case.\\

\begin{proposition}\label{conjunto I}
Let  $ \alpha $ and $ \beta $ be topological partial actions of $ G $ on $ X $ and $ Y$, respectively. If $ f, l: X \to Y $ are $G-$pmorphisms, then:
\begin{enumerate}
\item $ G $ acts partially and topologically on $ I = \{x \in X: f (x) = l (x) \} $.
\item The inclusion $ \sigma_I: I \to X $ with $ I = \{x \in X: f (x) = l (x) \} $ is an equalizer in the category $G-\mathbf{pActTop}$.
\end{enumerate}
\end{proposition}
\begin{proof}\hspace{0.5cm}
\begin{enumerate}
\item It is enough to prove that $ I $ is $ G-$pinvariant. In fact, if $ g \in G $ and $ y \in \alpha_g (I\cap X_{g^{-1}}) $, then $ y = \alpha_g(z) $ for some $z \in I\cap X_{g^{-1}}$. Thus,
\begin{align*}
f (y) = f (\alpha_g(z)) &= \beta_g(f (z)) = \beta_g(l (z)) = l (\alpha_g(z)) = l (y).
\end{align*}
Therefore, $ y \in I\cap X_{g}$ and so $ \alpha_g (I\cap X_{g^{-1}}) \subseteq I\cap X_{g}$.

\item It is evident that $ \sigma_I $ is a $G-$pmorphism such that $ f \circ \sigma_I = l \circ \sigma_I $. Now, suppose that $ G $ acts partially and topologically on a topological space $ J $ and $ k: J \to X $ is a $ G-$pmorphism such that $ f \circ k = l \circ k $. Then since the last equality implies that $k(J)\subseteq{I}$, we have that the function $ \lambda: J \to I $ defined as $ \lambda (j) = k (j) $ for each $ j \in J $ is well defined and is a $ G-$pmorphism such that $ \sigma_I \circ \lambda = k $. Finally, if $ \lambda ': J \to I $ is a $ G-$pmorphism such that $ \sigma_I \circ \lambda' = k $, then $ \lambda (j) = k (j) = (\sigma_I \circ \lambda ') (j) = \lambda' (j) $ for each  $ j \in J $. Thus, $ \lambda $ is unique and the following diagramm commutes:
$$
\xymatrix@=1.5cm{
I \ar[r]^{\sigma_I}  & X \ar@/_0.4pc/[r]_{l} \ar@/_{-0.4pc}/[r]^{f}&  Y\\
J  \ar@{-->}[u]^{\lambda} \ar[ur]_{k}& &
}
$$
\end{enumerate}
\end{proof}

\begin{theorem}
Let $ \alpha, \beta, \gamma $ be topological partial actions of $ G $ on $ X, Y, Z $ respectively. If $ f: X \to Z $ and $l: Y \to Z $ are $G-$pmorphisms, $\pi_X$ and $ \pi_Y $ are the projections  on $ X $ and $ Y $ respectively, $\sigma_P:P\to X\times Y$ is the equalizer of $f\circ \pi_X$ and $l\circ \pi_Y$, then the diagram
$$
\xymatrix@=1.5cm{
P \ar[r]^{\pi_X\circ\sigma_P} \ar[d]_{\pi_Y\circ\sigma_P} & X\ar[d]^{f}\\
Y\ar[r]_{l}& Z
}
$$
is a pullback and therefore, the category $G-\mathbf{pActTop}$ is finitely complete.
\end{theorem}
\begin{proof}[Proof]
To begin with, note that the diagram is commutative.
Now, $$ P = \{(x, y) \in X \times Y: f (x) = l (y) \}$$  and given $ (x, y) \in P $ we have
\begin{multline*}
(f \circ (\pi_X\circ\sigma_P)) (x, y) = f(\pi_X(\sigma_P(x, y))) =f (x)\\
= l (y) = l(\pi_Y(\sigma_P(x, y)))=(l \circ (\pi_Y\circ\sigma_P)) (x, y).
\end{multline*}
Also note that the topological partial action $\delta$ of $G$ of $P$ is the restriction of the topological partial action of $G$ of $X\times Y$ to $P$.
That is, $P_g=(X_g\times Y_g)\cap P$ and $\delta_g(x,y)=(\alpha_g(x),\beta_g(y))$ for each pair $g\in G$ and $(x,y)\in P_{g^{-1}}$.
Now, if
$$
\xymatrix@=1.5cm{
Q \ar[r]^{m} \ar[d]_{n} & X\ar[d]^{f}\\
Y\ar[r]_{l}& Z
}
$$
is a commutative diagram, let us define $ \rho: Q \to P $ by $ \rho (q) = (m (q), n (q)) $ for each $ q \in Q $.
Then, if $\varrho$ is a topological partial  action of $G$ on $Q$, for each pair $ g \in G $ and $ q \in Q $ we have
$$\rho(Q_g)=[m(Q_g)\times n(Q_g)]\cap P\subseteq(X_g\times Y_g)\cap P=P_g$$
and
\begin{align*}
\rho (\varrho_g(q)) &= (m (\varrho_g(q)), n (\varrho_g( q)))\\
&= (\alpha_g(m (q)), \beta_g(n (q)))\\
&= \delta_g(m (q),n (q))\\
&=\delta_g(\rho(q)).
\end{align*}

Also, $\rho$ is continuous since $ m $ and $ n $ are continuous. Thus $\rho$ is a $ G-$pmorphism.\\

Finally, $ (\pi_X\circ\sigma_P) \circ \rho = m $ and $ (\pi_Y\circ\sigma_P) \circ \rho = n $.
In addition, if $ \rho ': Q \to P $ is a $ G-$pmorphism such that $ (\pi_X\circ\sigma_P) \circ \rho' = m $ and $ (\pi_Y\circ\sigma_P) \circ \rho '= n $, then for each $ q \in Q $ we have
\begin{align*}
 \rho '(q) &= ((\pi_X\circ\sigma_P)(\rho' (q)), (\pi_Y\circ\sigma_P) (\rho '(q))) = (m (q), n (q) ) = \rho (q).
 \end{align*}
So $ \rho $ is unique.
$$
\xymatrix@=1.5cm{
Q\ar@{-->}[rd]^{\rho}\ar@/_{-1pc}/[rrd]^m \ar@/_{1pc}/[ddr]^n&&\\
&P \ar[r]^{\pi_X\circ\sigma_P} \ar[d]^{\pi_Y\circ\sigma_P} & X\ar[d]^{f}\\
&Y\ar[r]_{l}& Z
}
$$

Finally, the category $G-\mathbf{pActTop}$ is finitely complete by virtue of Theorem \ref{charprod}, Proposisition \ref{conjunto I}, and the existence of pullbacks given in this theorem.
\end{proof}


\subsection{Study of some colimits}

Recall that if $ \{X_i \}_{i \in I} $ is a collection of topological spaces, then the sum topology is considered on $ X = \bigsqcup_{i \in I} X_i $ (see \cite{munkres}).\\

\begin{proposition}\label{PartialActionSum}
Let $ \{X_i \}_{i \in I} $ be a collection of topological spaces and $ X = \bigsqcup_{i \in I} X_i $. If $ \alpha_i $ is a topological partial  action of $ G $ on $ X_i$ for every $i\in I$, then the function $ \alpha: \Gamma_\alpha \to X $ given by $ \alpha(g, x) = \alpha_{i_x}(g , x)$, for all 
$ (g, x) \in \Gamma_{\alpha}$, being $i_x\in I$ the (unique) index such that $x\in X_{i_x}$ and $\Gamma_{\alpha} = \{(h, y) \in G \times X: y \in X_{h^{-1}}\}$, is a topological  partial action of $ G $ on $X$. 
\end{proposition}
\begin{proof}
By assumption, for every $ i \in I $ there exists a class $ \{(X_{i, g}, \alpha_{i, g}) \}_{g \in G} $ such that for each $ g \in G $, $ X_ {i, g} \subseteq X_i$  is an open set, $ \alpha_{i, g}: X_{i, g^{- 1}} \to X_{i, g} $ is a homeomorphism, the set $ \Gamma_{\alpha_i} = \{(g, x) \in G \times X_i: x \in X_{i, g^{- 1}}\}$ is open in $ G \times X_i $ and the function $ \alpha_i: \Gamma_{\alpha_i} \to X $ given by $ \alpha_i (g, x) = \alpha_{i, g} (x) $ is continuous and also a partial action  of $ G $ on $ X_i $.
Now, for each $ g \in G $, let us consider $ X_{g} = \bigsqcup_ {i \in I} X_{i, g} $.\\

So, if $ x \in X_ {g^{- 1}} $, there exists a unique $ i_x \in I $ such that $ x \in X_ {i_x, g^{- 1}} $.
Therefore, the function $ \alpha_ {g}: X_{g^{- 1}} \to X_{g} \subseteq X $ given by $ \alpha_ {g} (x) = \alpha_{i_x, g} (x ) $ is well defined and also a bijection.
Now, if $ x, y \in X_{g ^ {- 1}} $ with $ x = y $, then $ i_x = i_y $ and therefore $ \alpha_{g} (x) = \alpha_{i_x, g} (x) = \alpha_{i_y, g} (y) = \alpha_ {g} (y) $.
On the other hand, if $ \alpha_{g} (x) = \alpha_ {g} (y) $, then $ x, y \in X_{i_0, g^{- 1}} $ for a unique $ i_0 \in I $.
So $ \alpha_{i_0, g} (x) = \alpha_ {i_0, g} (y) $ where $ x = y $.
The case $ x \in X_ {i, g^{- 1}} $ and $ y \in X_ {j, g^{- 1}} $ with $ i \neq j $ can not occur, because $ X_i $ and $ X_j $ are disjoint when $ i \neq j $.
Hence $ \alpha_g $ is injective for each $ g \in G $.
Also, given $ y \in X_ {g} $, we have $ y \in X_ {i_y, g} $ for a unique $ i_y \in I $, and therefore,  we may define $ x = \alpha_ {i_y, g}^{-1}(y) \in X_{i_y, g^{- 1}} $.
Thus, $ \alpha_g (x) = \alpha_ {i_y, g} (\alpha_{i_y, g}^{- 1}(y)) = y $ and so $ \alpha_g $ is surjective for each $ g \in G $.
Since $ X_ {i, g} $ is open for every $ i \in I $, then $ X_g $ is an open set for each $ g \in G $ (with the sum topology).
Also, since $ \alpha_ {i, g} $ is a homeomorphism for every $ i \in I $, it follows that $ \alpha_g $ is a homeomorphism for each $ g \in G $.\\

Since items 1, 2 and 3 of Definition \ref{defpartop} are true (see \cite{Avila1} for details), only remains to prove  that the set $\Gamma_{\alpha}$ is open and that the corresponding function $ \alpha: \Gamma_{\alpha} \to X $ is well defined and is continuous.\\

If $ (g, x) \in \Gamma_{\alpha}$, then $ x \in X_{g^{- 1}} $. Therefore, $ x \in X_{i_x, g^{- 1}} $ for a unique $ i_x \in I $ and thus $ (g, x) \in \Gamma_{\alpha_ {i_x}} $. Since $\Gamma_{\alpha_ {i_x}}$  is an open set in $ G \times X_ {i_x} $, there are open neighborhoods $M\subseteq G$ and $N\subseteq X_{i_x}$ of $g$ and $x$ respectively such that $M\times N\subseteq\Gamma_{\alpha_ {i_x}}$. Let us notice that $\Gamma_{\alpha_{i_x}} \subseteq \Gamma_\alpha $  because  $X_{i_x} \subseteq X$.  Consequently,  $M\times N \subseteq \Gamma_\alpha $ and therefore, $ \Gamma_{\alpha} $ is an open set in $ G \times X $.\\

Suppose that $ (g, x) \in \Gamma_ \alpha $ and $V\subseteq X$ is an open neighborhood of $\alpha(g,x)$. Then $ x \in X_ {g^{-1}} $ and thus $ x \in X_ {i_x, g^{-1}} $ for a unique $ i_x \in I $. Hence, $ (g, x) \in \Gamma_{\alpha_{i_x}} $ and 
$\alpha (g, x) = \alpha_{i_x}(g, x) =\alpha_{i_x, g}(x) =\alpha_ {g} (x) $
is well defined. Since $ \alpha_{{i_x}}: \Gamma_{\alpha_i} \to X $ is continuous in  $ (g, x) \in \Gamma_{\alpha_i} $, then there exists $U_{{i_x}}$ open in $\Gamma_{\alpha_{i_x}}$ such that  $(g,x)\in U_{{i_x}}$ and $\alpha_{i_x}(U_{{i_x}})\subseteq V$.  Since $\Gamma_\alpha$ is an open set containing $\Gamma_{\alpha_{i_x}}$, then  $U_{{i_x}}$ is  open in $\Gamma_{\alpha}$ and
$\alpha(U_{{i_x}})=\alpha_{i_x}(U_{{i_x}})\subseteq V$. So $ \alpha: \Gamma_{\alpha} \to X $ is continuous.
\end{proof}

The following result  characterizes the co-products in $G-\mathbf{pActTop}$.\\

\begin{theorem}
Let $ \{X_i \}_{i \in I} $ be a collection of topological spaces. Then the topological  partial action given by Proposition \ref{PartialActionSum} on 
  $ X = \bigsqcup_{i \in I} X_i $ characterizes the co-product of the collection $ \{X_i \}_{i \in I} $ (with the sum topology) in the category $G-\mathbf{pActTop}$.
\end{theorem}
\begin{proof}
With the same notation of Proposition \ref{PartialActionSum}. If $ \sigma_i: X_i \to X $ is the inclusion of $ X_i $ into $ X $ given by $ \sigma_i (x) = x $, then $ \sigma_i $ is continuous, $\sigma_i(X_{i,g})=X_{i,g}\subseteq X_g$ and $ \sigma_i (\alpha_g (x)) = \alpha_{i,g}(x) = \alpha_{i,g}(\sigma_i (x)) $ for each pair $ g \in G $ and $ x \in X_{i,g^{-1}} $. Therefore,  $ \sigma_i $ is a $G-$pmorphism.\\

If  $ \beta$ is a topological partial action of $G$ on $Y$ and $ f_i: X_i \to Y $ is a $G-$pmorphism for every $ i \in I $, then $ f_i $ is continuous, $f_i(X_{i,g})\subseteq Y_g$ and $f_i(\alpha_{i,g}(x))=\beta_g(f_i (x))$ for each pair $ g \in G $ and $ x \in X_{g^{-1}} $. Now, for each $ x \in X $ there exists a unique $ i_x \in I $ such that $ x \in X_{i_x} $ and $ f_{i_x} (x) \in Y $. Thus,  the function $ f: X \to Y $ given by $ f (x) = f_{i_x} (x) $ is well defined and continuous. Indeed, suppose that $x\in X$ and $U\subseteq Y$ is an open neighborhood of $f(x)$. Then, there exists a unique $ i_x \in I $ such that $ x \in X_{i_x} $ and $f(x)=f_{i_x}(x)$. Since $f_{i_x}:X_{i_x} \to Y $ is continuous, then there exists $V_{i_x}$ open in $X_{i_x}$ such that $x\in V_{i_x}$ and $f_{i_x}(V_{i_x})\subseteq U$. Let us notice that $V_{i_x}$ is open in $X$ and $f(V_{i_x})=f_{i_x}(V_{i_x})\subseteq U$. In addition,
$$f(X_g)=f\biggl(\bigsqcup_{i\in I}X_{i,g}\biggr)=\bigcup_{i\in I}f_i(X_{i,g})\subseteq\bigcup_{i\in I}Y_g=Y_g$$
and
$$ f (\alpha_g(x)) = f_{i_x} (\alpha_{i_x,g}(x)) = \beta_g(f_{i_x} (x)) = \beta_g(f (x)) $$ for each pair $ g \in G $ and $ x \in X_{g^{-1}} $. That is, $ f $ is a $G-$pmorphism.\\

Finally, $ f $ is the unique function such that $ f \circ \sigma_i = f_i $ for every $ i \in I $. Indeed, $ (f \circ \sigma_i) (x) = f (\sigma_i (x)) = f (x) = f_i (x) $ for each $ x \in X_i $. That is, the diagram below is commutative:
$$
\xymatrix@=1.5cm{
\bigsqcup_{i \in I}X_i\ar@{-->}[dr]^f& \\
X_i \ar[u]^{\sigma_i} \ar[r]_{f_i}&Y
}
$$
Also, if $ h: X \to Y $ is a $G-$pmorphism such that $ h \circ \sigma_i = f_i $ for  every  $ i \in I $, then $ f (x) = f_{i_x} (x) = h (\sigma_ {i_x} (x)) = h (x) $ for each $ x \in X $. Hence, $ h = f $.
\end{proof}

To study of co-equalizers and pushouts, it is necessary  to consider the following notion.\\

\begin{definition}\label{gpcong}
Given a topological partial action $\alpha$ of $G$ on $X$, an equivalence relation $\varepsilon_X$ on $X$, with equivalence classes denoted by $\overline{x}$, for each $x \in X$, is said to be a $G-$pcongruence if satisfies the following conditions:
\begin{enumerate}
  \item For every $x \in X$ and each open set $U \subseteq X$ such that $x \in U$, it holds that $\overline{x} \subseteq U$.
  \item If $x, y \in X_{g^{-1}}$ such that $\overline{x} = \overline{y}$, then $\overline{\alpha_g(x)} = \overline{\alpha_g(y)}$.
\end{enumerate}
\end{definition}

\begin{proposition}Using the above notation, let $\varepsilon_X$ be a $G-$pcongruence on $X$ and $\mathcal{Q}$ the quotient topology on $X/\varepsilon_X$. Then, the quotient function  $\tau_X: X \to X/\varepsilon_X$ is open and continuous.\\
\end{proposition}

Next, we will see how every $G-$pcongruence  $\varepsilon_X$ induces a topological partial action on $X/\varepsilon_X$. From now on, we will consider the quotient topology $\mathcal{Q}$ on $ X / \varepsilon_X$. So, the quotient function $\tau_X: X \to X / \varepsilon_X$ is continuous and open.\\

\begin{proposition}\label{quotient}
Let $ \alpha $ be a topological partial action of $ G $ on $ X $. If $ \varepsilon_X $ is a $G-$pcongruence on $ X $, then $ G $ acts partially and topologically on $ X /\varepsilon_X $.
\end{proposition}
\begin{proof}
With the previous notation, let us consider the class $ \{(\Sigma_g, \Lambda_g) \}_{g \in G} $ where $ \Sigma_g = \tau_X(X_g) $ and $ \Lambda_g: \Sigma_{g^{-1}} \to \Sigma_{g} $ is defined by $ \Lambda_g (\overline {x}) = \overline {\alpha_g (x)} $ for each $\overline{x}\in \Sigma_{g^{-1}}$.
\begin{enumerate}
\item Note that if $ g \in G $, then $X_g\subseteq X$ is open, and therefore, $\Sigma_g = \tau_X (X_g) \in \mathcal{Q}$ is an open subset of $ X / \varepsilon_X$.\\

Additionally, if $ g \in G $ and $ \overline{x}, \overline{y}\in \Sigma_{g^{-1}} $ are such that $ \overline {x} = \overline {y} $, then $ \overline{\alpha_g (x)} = \overline {\alpha_g (y)} $, and thus, $ \Lambda_g (\overline {x}) = \Lambda_g (\overline {y}) $, which means $ \Lambda_g $ is well-defined. Conversely, $ \Lambda_g $ is injective, as
$ \Lambda_g (\overline {x})=\Lambda_g (\overline {y})$ implies $\overline{\alpha_g (x)} = \overline {\alpha_g (y)}$ and so, $\overline{x}=\overline{\alpha_{g^{-1}} (\alpha_g (x))} = \overline {\alpha_{g^{-1}}(\alpha_g (y))}=\overline{y}$. Moreover, if $ g \in G $ and $ \overline{y} \in\Sigma_g $, then $\tau_X(y)\in \tau_X(X_g)$, and therefore,  there exists $z \in X_g $ such that $\tau_X(y)=\tau_X(z)$. Since $ \alpha_g $ is surjective, there exists $ x \in X_{g^{-1}} $ such that $ z = \alpha_g (x) $, and consequently, $ \overline {x} \in \Sigma_{g^{-1}} $. Hence, $ \overline{y} = \overline{z}= \overline{\alpha_g (x)}=\Lambda_g (\overline{x})$, and so, $ \Lambda_g $ is surjective.\\

Finally, for $ g \in G $, $ \overline{x} \in \Sigma_{g^{-1}} $, and $ U \subseteq \Sigma_{g} $ an open set such that $  \Lambda_g (\overline {x}) \in U $, we have $\tau_X(\alpha_g (x)) \in U $. Since $ \tau_X $ restricted to $X_{g}$ is continuous, there is an open $V\subseteq X_{g}$ such that $\alpha_g (x)\in V$ and $\tau_X(V) \subseteq U$. Since $\alpha_g$ is continuous, there exists an open $ W\subseteq X_{g^{-1}}$ such that $ x \in W $ and $ \alpha_g (W) \subseteq V $. Thus, $ \tau_X (W)\in \mathcal{Q}$ is open with $ \overline {x} \in \tau_X (W) $ and $ \Lambda_g (\tau_X (W))= \tau_X (\alpha_g(W)) \subseteq  \tau_X (V) \subseteq U $. Therefore, $ \Lambda_g $ is continuous, and since $ \Lambda_g^{-1} = \Lambda_{g{^{-1}}}$ is also continuous, $ \Lambda_g $ is a homeomorphism for each $ g \in G $.

\item $ \Sigma_e = \tau_X (X_e) = \tau_X (X) = X / \varepsilon_X $, and for each $ \overline {x} \in X / \varepsilon_X $, we have
$ \Lambda_e (\overline {x}) = \overline {\alpha_e (x)} = \overline {x}$, meaning $ \Lambda_e = i_{X / \varepsilon_X} $.

\item If $ g, h \in G $ and $ x \in X $ are such that $ \overline {x} \in \Lambda_h^{-1} (\Sigma_h \cap \Sigma_{g^{-1}}) $, then $$ \tau_X(\alpha_h (x))=\Lambda_h (\overline {x}) \in \Sigma_h \cap \Sigma_{g^{-1}} = \tau_X(X_h) \cap  \tau_X(X_{g^{-1}}).$$ 
So, there exist $z\in X_h$ and $z'\in X_{g^{-1}}$ such that $\overline{z}=\overline{\alpha_h(x)}$ and $\overline{z'}=\overline{\alpha_h(x)}$. 
By 2. from Definition \ref{gpcong}, we have that $\overline{z}\subseteq X_h\cap X_{g^{-1}}$, and thereby, $z\in \alpha_h(X_{(gh)^{-1}})$. That is, 
$\Lambda_h (\overline {x})=\overline{z}\in \tau_X(\alpha_h(X_{(gh)^{-1}}))=\Lambda_h(\tau_X(X_{(gh)^{-1}}))$ or equivalently, 
$\overline {x}\in \tau_X(X_{(gh)^{-1}})=\Sigma_{(gh)^{-1}}$.
\item Now, for $ g, h \in G $ and $ x \in X $ such that  $ \overline{x} \in \Lambda_h^{- 1} (\Sigma_h \cap \Sigma_{g^{ -1}})$,
$$ (\Lambda_g \circ \Lambda_h) (\overline{x}) = \Lambda_g(\Lambda_h (\overline{x})) =  \Lambda_g (\overline{\alpha_h (x)}) = \overline{\alpha_g (\alpha_h (x))} = \overline{\alpha_{gh} (x)} = \Lambda_{gh} (\overline{x}). $$

\item Let us suppose that $ (h, \overline{y}) \in \Gamma_{\Lambda} = \{(g,\overline{x}) \in G \times (X / \varepsilon_X): \overline{x} \in \Sigma_{g^{-1}} \} $. Since $\tau_X(y)=\overline{y} \in \Sigma_{h^{-1}}=\tau_X(X_{h^{-1}})$, there exists $z\in X_{h^{-1}}$ such that $\tau_X(y)=\tau_X(z)$. Consequently, $(h, z) \in \Gamma_\alpha$, which is open, and thereby, there exist open $ M \subseteq G $ and $ N \subseteq X $ such that $ (h, z) \in M\times N \subseteq \Gamma_\alpha $. Moreover, if $(t,\overline{w})\in M \times \tau_X (N)$, then there exists $w'\in N$ such that $\tau_X(w)=\tau_X(w')$. Hence,  $(t,w')\in M\times N \subseteq \Gamma_{\alpha}$, and  
therefore, $w'\in X_{t^{-1}}$. Thus, $\overline{w}=\overline{w'}\in \tau_X(X_{t^{-1}})=\Sigma_{t^{-1}}$. That is, $(t,\overline{w})\in \Gamma_{\Lambda}$, and thereby, $ M \times \tau_X (N) \subseteq \Gamma_{\Lambda}$. Hence, $\Gamma_\Lambda$ is open.
\item The function $ \Lambda: \Gamma_\Lambda \to X / \varepsilon_X $ defined by $ \Lambda (g, \overline {x}) = \Lambda_g (\overline {x}) $ is continuous as $ \Lambda (g, \overline {x}) = (\tau_X \circ \alpha) (g, x) $ and  $ \tau_X $ is open and continuous, and $ \alpha $ is continuous.
\end{enumerate}
Consequently, $ \Lambda $ is a topological partial action of $ G $ on $ X / \varepsilon_X $.
\end{proof}

\begin{definition}[\cite{Adamek}]
Given an object $A$ in a category $\mathbf{C}$, it is said that an object $B$ in $\mathbf{C}$ is a quotient object of $A$ if there exists an epimorphism $f: A \to B $. \\  
\end{definition}
\begin{theorem}
The topological partial action $\Lambda$ of $G$ on $X/\varepsilon_X$ from Proposition \ref{quotient} is a quotient object of the topological partial action $\alpha$ of $G$ on $X$.
\end{theorem}
\begin{proof}
Indeed,  $\tau_X$ is continuous with the quotient topology  over $X / \varepsilon_X$. Now, for each pair $ g \in G $ and $ x \in X_{g^{-1}}$, we have that $\tau_X (X_g)=\Sigma_g $  and
$$ \tau_X (\alpha_g( x)) = \overline {\alpha_g( x)} = \Lambda_g(\overline {x}) = \Lambda_g(\tau_X (x)).$$
In addition, if $h,l: X / \varepsilon_X \to Y$ are $G-$pmorphisms such that
$h \circ \tau_X = l \circ \tau_X$, then 
$$h(\overline{x})=h(\tau_X(x))=l(\tau_X(x))=l(\overline{x}).$$
$$
\xymatrix@=1.5cm{
X \ar[r]^{\tau_X} & x/\varepsilon_X \ar@/_0.4pc/[r]_{l} \ar@/_{-0.4pc}/[r]^{h} &  Y
}
$$
Thereby, 
the function $\tau_X:X \to X/\varepsilon_X$ is an epimorphism in the category $G-$\textbf{pActTop}, and thus, $X/\varepsilon_X$ is a quotient object of $X$.
\end{proof}

Let us notice that if $ f: X \to Y $ is a $G-$pmorphism, then the condition $ f (x) \in Y_{g^{- 1}} $ for some pair $ g \in G $ and $ x \in X $ does not imply that $ x \in X_{g^{-1}} $. This situation is not considered in topological global actions since it is always true that $ X_g = X $ and $ Y_g = Y $ for each $ g \in G $. However, in the case of topological partial actions, this distinction is important in the construction of co-equalizers and pushouts.\\

\begin{definition}
Let $ \alpha $ and $ \beta $ be topological partial actions of $ G $ on $ X $ and $ Y $ respectively. We say that a $G-$pmorphism $ f: X \to Y $ is strong if $ f (x) \in Y_g $ implies that $ x \in X_g $ for each $ g \in G $.\\
\end{definition}

\begin{proposition}
Let $ \alpha $ be a topological partial action of $ G $ on $ X $. If $ \varepsilon_X $ is a $G-$pcongruence on $ X $ and $\Lambda$ is the topological partial action of $ G $ on $ X/\varepsilon_X $ from Proposition \ref{quotient}, then the quotient function $ \tau_X: X \to X / \varepsilon_X $ is a strong $G-$pmorphism.
\end{proposition}

\begin{proof}
If $g\in G$ and $\tau_X(x)\in \Lambda_g$, then $\tau_X(x)\in \tau_X(X_g)$. So, there exists $y\in X_g$ such that $\tau_X(x) = \tau_X(y)$. By 1. from Definition \ref{gpcong}, we have $x \in \overline{x} = \overline{y} \subseteq X_g$.
\end{proof}

Observe that if $ \alpha $ and $ \beta $ are topological partial actions of $ G $ on $ X $ and $ Y $ respectively, then the strong $G-$pmorphisms $ f, l: X \to Y $ induce a relation $ R_Y $ on $ Y $ given by:\\

\emph{For $y_1, y_2 \in Y$, we say that $(y_1, y_2) \in R_Y$ if and only if there exists $g \in G$ such that $y_1, y_2 \in Y_g$ and there exists $ x \in X_{g} $ such that $y_1 = f(x)$, $y_2 = l(x)$, and for all open $V \subseteq X_g$ containing $x$, we have that $f(V) = l(V)$.}\\

If $\Delta_Y = \{(y,y) : y \in Y\}$ and $\Omega_Y \subseteq Y \times Y$ is defined by:
\begin{center}
\emph{$(u, v) \in \Omega_Y$ if there exist $y_1, \dots, y_n \in Y$ such that $u = y_1$, $v = y_n$ and $(y_1, y_2), \dots, (y_{n-1}, y_n) \in R_Y \cup R_Y^{-1} \cup \Delta_Y$,}
\end{center}
then $\varepsilon_Y = R_Y \cup R_Y^{-1} \cup \Delta_Y \cup \Omega_Y$
is the minimal equivalence relation containing the relation $ R_Y $ induced by $f$ and $l$. This leads to the following result.\\

\begin{proposition}
Let $ \alpha $ and $ \beta $ be topological partial actions of $ G $ on $ X $ and $ Y $ respectively. If $ f, l: X \to Y $ are strong $G-$pmorphisms and $ \varepsilon_Y$ is the above equivalence relation on $Y$ induced by $f$ and $l$, then $\varepsilon_Y $ is a $ G-$pcongruence on $ Y $.
\end{proposition}
\begin{proof}\hspace{0.5cm}
\begin{enumerate}
    \item Let $y \in Y$ and an open $U \subseteq Y$ such that $y \in U$. If $z \in \overline{y}$, then:
    \begin{enumerate}
        \item If $(z, y) \in R_Y$, there exists $g \in G$ such that $z, y \in Y_g$ and there exists $ x \in X_g $ such that $z = f(x)$, $y = l(x)$, and for every open $V \subseteq X_g$ with $x \in V$, we have $f(V) = l(V)$. Consequently, $l(x) \in U$ and since $l$ is continuous, there exists an open $W_x \subseteq X$ with $x \in W_x$ and $l(W_x) \subseteq U$. Taking $V = W_x \cap X_g$, it follows that $z = f(x) \in f(W_x \cap X_g) = l(W_x \cap X_g) \subseteq l(W_x) \subseteq U$.
        \item If $(z, y) \in R_Y^{-1}$, then $(y, z) \in R_Y$ and the proof is analogous to the previous case.
        \item If $(z, y) \in \Delta_Y$, then $z = y \in U$.
        \item If $(z, y) \in \Omega_Y$, then there exists $y' \in Y$ that satisfies one of the following cases:
        \begin{enumerate}
            \item If $(z, y'), (y', y) \in R_Y$, there exist $g, h \in G$ such that $z, y' \in Y_g$ and $y, y' \in Y_h$ and $x \in X_g$ and $t \in X_h$ such that $z = f(x)$, $y = l(t)$, and $y' = f(t) = l(x)$. Thereby, $l(t) \in U$, and since $l$ is continuous, there exists $W_t \subseteq X$ open such that $t \in W_t$ and $l(W_t) \subseteq U$. Since $(y', y) \in R_Y$, then $l(x) = f(t) \in f(W_t \cap X_h) = l(W_t \cap X_h) \subseteq l(W_t) \subseteq U.$ Similarly, there exists $W_x \subseteq X$ open with $x \in W_x$ and $l(W_x) \subseteq U$, and since $(z, y') \in R_Y$, it follows that $z = f(x) \in f(W_x \cap X_g) = l(W_x \cap X_g) \subseteq l(W_x) \subseteq U.$
            \item If $(z, y') \in R_Y$ and $(y', y) \in R_Y^{-1}$, then $(y, y') \in R_Y$ and the argument follows similarly from the previous case.
            \item If $(z, y') \in R_Y$ and $(y', y) \in \Delta_Y$, then $y = y'$ and there exists $g \in G$ such that $z, y' \in Y_g$. Furthermore, there exists $x \in X_g$ with $z = f(x)$, $y' = l(x)$. Thus, $l(x) \in U$. Since $l$ is continuous, there exists $W_x \subseteq X$ open with $x \in W_x$ and $l(W_x) \subseteq U$. So,  $z = f(x) \in f(W_x \cap X_g) = l(W_x \cap X_g) \subseteq l(W_x) \subseteq U$ because $(z, y') \in R_Y$.
            \item If $(z, y'), (y', y) \in R_Y^{-1}$, then $(y, y'), (y', z) \in R_Y$ and the proof follows similarly as in case i.
            \item If $(z, y') \in R_Y^{-1}$ and $(y', y) \in \Delta_Y$, then $y = y'$ and $(y', z) \in R_Y$, and the proof is analogous to case iii.
            \item If $(z, y'), (y', y) \in \Delta_Y$, then $z = y' = y \in U$.
            \item If $(z, y') \in R_Y^{-1}$ and $(y', y) \in R_Y$, then it is similar to the case i.
            \item If $(z, y') \in \Delta_Y$ and $(y', y) \in  R_Y$, then it is similar to the case iii.
            \item If $(z, y') \in \Delta_Y$ and $(y', y) \in R_Y^{-1}$, then it is similar to the case iii.
        \end{enumerate}
    \end{enumerate}
\item Given $g \in G$ and $y_1, y_2 \in Y_{g^{-1}}$ such that $\overline{y_1} = \overline{y_2}$, we have the following cases:
      \begin{enumerate}
        \item If $ (y_1 , y_2 ) \in R_Y$, then there exists $ x \in X_{g^{-1}}$ such that $ y_1 = f(x) $ and $ y_2 = l(x) $. Hence, we have $ \beta_g( y_1) = \beta_g(f(x)) = f(\alpha_g( x)) $ and $ \beta_g( y_2) = \beta_g( l(x) ) = l(\alpha_g( x))$. Now, if $V \subseteq X_g$ is an open set such that $\alpha_g(x) \in V$, then $x \in \alpha_{g^{-1}}(V)$, where $\alpha_{g^{-1}}(V) = \alpha_g^{-1}(V)$ is an open set of $X_{g^{-1}}$, since $\alpha_{g^{-1}} : X_{g} \to X_{g^{-1}}$ is a homeomorphism. So, 
$$\beta_{g^{-1}}(f(V)) = f(\alpha_{g^{-1}}(V)) = l(\alpha_{g^{-1}}(V)) = \beta_{g^{-1}}(l(V)).$$
            Since $\beta_{g^{-1}}(f(V)), \beta_{g^{-1}}(l(V)) \subseteq Y_{g^{-1}}$, we apply $\beta_{g}$ and obtain that $f(V) = l(V)$. Therefore,  $(\beta_g(y_1), \beta_g(y_2)) \in R_Y \subseteq \varepsilon_Y$, and consequently, $\overline{\beta_g(y_1)} = \overline{\beta_g( y_2)}.$
        \item If $ (y_1 , y_2 ) \in R_Y^{-1}$, then $ (y_2 , y_1 ) \in R_Y$. By analogy to the previous case, $(\beta_g(y_2), \beta_g(y_1)) \in R_Y$, and thereby, $(\beta_g(y_2), \beta_g(y_1)) \in R_Y^{-1} \subseteq \varepsilon_Y$. Hence,  $\overline{\beta_g(y_1)} = \overline{\beta_g( y_2)}.$
        \item If $ (y_1 , y_2 ) \in \Delta_Y$, then $y_1 = y_2$, and therefore, $\beta_g(y_1) = \beta_g(y_2)$. So, $(\beta_g(y_1), \beta_g(y_2)) \in \Delta_Y \subseteq \varepsilon_Y$, and thus, $\overline{\beta_g(y_1)} = \overline{\beta_g( y_2)}.$
        \item If $ (y_1 , y_2 ) \in \Omega_Y$, then there exists at least one $y' \in Y$ that satisfies one of the following cases:
  \begin{enumerate}
\item If $(y_1, y'), (y', y_2) \in R_Y$, then by case $a)$,
$$\overline{\beta_g(y_1)} = \overline{\beta_g(y')} = \overline{\beta_g(y_2)}.$$
\item If $(y_1, y') \in R_Y$ and $(y', y_2) \in R_Y^{-1}$, it follows by cases $a)$ and $b)$.
\item If $(y_1, y') \in R_Y$ and $(y', y_2) \in \Delta_Y$, it follows by cases $a)$ and $c)$.
\item If $(y_1, y'), (y', y_2) \in R_Y^{-1}$, it follows by case $b)$.
\item If $(y_1, y') \in R_Y^{-1}$ and $(y', y_2) \in \Delta_Y$, it follows by cases $b)$ and $c)$.
\item If $(y_1, y'), (y', y_2) \in \Delta_Y$, it follows by case $c)$.
\item If $(y_1, y') \in R_Y^{-1}$ and $(y', y_2) \in  R_Y$, then it is similar to the case ii.
\item If $(y_1, y') \in \Delta_Y$ and $(y', y_2) \in  R_Y$, then it is similar to the case iii.
\item If $(y_1, y') \in \Delta_Y$ and $(y', y_2) \in  R_Y^{-1}$, then it is similar to the case v.
\end{enumerate}
\end{enumerate}
\end{enumerate}
\end{proof}
Before presenting the following result, we define the category $G-\mathbf{spActTop}$ in which the objects are all ordered pairs $(X,\alpha)$ where $X$ is a topological space and $\alpha$ is a  topological partial action of $G$ on  $X$, and the morphisms are all strong $G-$pmorpisms. Note that $G-\mathbf{spActTop}$ has co-products since the inclusions are strong $G-$pmorpisms.\\
\begin{theorem}
Let $ \alpha $ and $ \beta $ be topological partial actions of $ G $ on  $ X $ and $ Y $ respectively. If $ f, g: X \to Y $ are strong $G-$pmorphisms and $ \varepsilon_Y $ is the minimal equivalence relation induced by $ f $ and $ g $, then $ \tau_Y: Y \to Y / \varepsilon_Y $ is a co-equalizer of $ f $ and $ g $ in the category $G-\mathbf{spActTop}$.
\end{theorem}
\begin{proof}
First note that $ (\tau_Y \circ f) (x) = \overline {f (x)} = \overline {g (x)} = (\tau_Y \circ g) (x) $ for each $ x \in X $. Now, if $ \{(\Sigma_t, \Lambda_t) \}_{t \in G} $ and $ \{(Z_t, \gamma_t) \}_{t \in G} $ are the classes  that define the topological partial actions of $G$ on  $ Y/\varepsilon_Y $  and  $ Z $ respectively and $ h: Y \to Z $ is a strong $ G-$pmorphism such that $ h \circ f = h \circ g $, let us consider the function $ \eta: (Y / \varepsilon_Y) \to Z $ as $ \eta (\overline {y}) = h (y) $ for each $ \overline {y} \in Y / \varepsilon_Y $. If $ \overline{y_1} = \overline{y_2} $, then there exists $ x \in X $ such that $ y_1 = f (x) $ and $ y_2 = g (x) $ where $$ \eta (\overline{y_1}) = h (f (x)) = h (g (x)) = \eta (\overline{y_2} ).$$
Consequently, $ \eta $ is well defined. Also, $\eta(\Sigma_t)=\eta(\tau_Y(Y_t))=h(Y_t)\subseteq Z_t$ and
\begin{align*}
\eta (\Lambda_t(\overline{y})) &=\eta (\Lambda_t(\tau_Y(y)))\\
&=\eta (\tau_Y(\beta_t( y))) \\
&=\eta (\overline {\beta_t( y)})\\
&= h (\beta_t( y))\\
&= \gamma_t( h (y)) \\
&= \gamma_t(\eta (\overline {y})),
\end{align*}
for each pair $ t \in G $ and $ \overline {y} \in \Sigma_{t^{-1}} $. Moreover, we have that $ \eta $ is continuous since $ h $ is continuous and $ \tau_Y $ is open.  In addition, if $g\in G$ and $\eta(\overline{y})\in Z_g$, then $h(y)\in Z_g$, and thereby, $y\in Y_{g}$. So, $\overline{y}=\tau_Y(y)\in \tau_Y(Y_g)= \Sigma_g$, and therefore, $ \eta $ is a strong $ G-$pmorphism.\\

Finally, $ (\eta \circ \tau_Y) (y) = \eta (\overline {y}) = h (y) $ for each $y\in Y $ and if $ \eta':( Y / \varepsilon_Y ) \to Z $ is a strong $ G-$pmorphism such that $ \eta '\circ \tau_Y = h $, then $$ \eta' (\overline {y}) = (\eta'(\tau_Y (y)) = h (y) = \eta (\overline{y}).$$
Hence, $ \eta $ is unique.
$$
\xymatrix@=1.5cm{
X \ar@/_0.4pc/[r]_{g} \ar@/_{-0.4pc}/[r]^{f}& Y \ar[r]^{\tau_Y} \ar[dr]_{h} &  Y/\varepsilon_Y \ar@{-->}[d]^{\eta}\\
&   & Z
}
$$
\end{proof}

\begin{theorem}
Let $ \alpha, \beta, \gamma $ be topological partial actions of $ G $ on $ X, Y, Z $ respectively. If $ f: X \to Y $ and $ g: X \to Z $ are a strong $G-$pmorphism, $ W = Y \sqcup Z $, $ \sigma_Y: Y \to W $ and $ \sigma_Z: Z \to W $ are the inclusions of $ Y $ and $ Z $ in $W$ respectively, then $ \sigma_Y \circ f $ and $ \sigma_Z \circ g $ induce a $ G-$pcongruence $ \varepsilon_W $ on $ W $, and therefore, the diagram
$$
\xymatrix@=1.5cm{
X \ar[r]^{f} \ar[d]_{g} & Y\ar[d]^{\tau_W\circ \sigma_Y}\\
Z\ar[r]_{\tau_W\circ \sigma_Z}& W/\varepsilon_W
}
$$
is a pushout. Thus, category $G-\mathbf{spActTop}$ is co-complete.
\end{theorem}
\begin{proof}
First note that the diagram is commutative. Indeed, given $ x\in X $ we have
$$ [(\tau_W \circ \sigma_Y) \circ f] (x) = \overline{f (x)} = \overline{g (x)} = [(\tau_W \circ \sigma_Z) \circ g] (x) $$
Now, if
$$
\xymatrix@=1.5cm{
X \ar[r]^{f} \ar[d]_{g} & Y\ar[d]^{t}\\
Z\ar[r]_{s}& R
}
$$
is a commutative diagram. Let us consider  the function $ \varsigma : (W / \varepsilon_W) \to R $  as
$$\varsigma(\overline{w})=
\left\{
\begin{array}{ll}
  t(w), & \text{if} \hspace{0.2cm} w\in Y\\\\
    s(w), & \text{if} \hspace{0.2cm} w\in Z.
     \end{array}
\right.
$$
Additionally, $ \varsigma $ is well defined and is continuous. Also, if
$ \{(W_h,\omega_h)\}_{h\in G} $, $ \{(\Omega_h,\Delta_h)\}_{h\in G} $ and $\{(R_h,\chi_h)\}_{h\in G}  $
are the classes that define the topological partial  actions of $ G $ on
$ W $, $ W / \varepsilon_W $ and $ R $,
respectively, then for each $ h \in G $ we have that
$$\varsigma(\Omega_h)=\varsigma(\tau_W(Y_h\sqcup Z_h))=t(Y_h)\cup s(Z_h)\subseteq R_h\cup R_h=R_h$$
and
\begin{center}
$\varsigma(\Delta_h(\overline{w}))=\varsigma(\overline{\omega_h(w)})=t(\omega_h(w))=\chi_h(t(w))=\chi_h(\varsigma(\overline{w})),$ if $w\in Y$,
\end{center}
\begin{center}
$\varsigma(\Delta_h(\overline{w}))=\varsigma(\overline{\omega_h(w)})=s(\omega_h(w))=\chi_h(s(w))=\chi_h(\varsigma(\overline{w})),$ if $w\in Z$.
\end{center}
Now, suppose $h\in G$ and $\varsigma (\overline{w})\in R_h$. If $w\in Y$, then $t(w)\in R_h$, and thus, $w\in Y_g$ because $t$ is a strong $G-$pmorphism. Consequently, $\overline{w}=(\tau_W\circ \sigma_Y)(y)\in \tau_W(Y_h)=\Omega_h$. If $w\in Z$, then it is similar since $s$ is also a strong $G-$pmorphism. Thereby, $ \varsigma $ is a strong $ G-$pmorphism.\\

On the other hand, for each pair $ y \in Y $ and $ z \in Z $ we have
\begin{center}
$[\varsigma\circ(\tau_W\circ \sigma_Y)](y)=\varsigma(\overline{y})=t(y)$ and $[\varsigma\circ(\tau_W\circ \sigma_Z)](z)=\varsigma(\overline{z})=s(z)$.
\end{center}
Finally, if $ \varsigma ':( W / \varepsilon_W) \to R $ is a strong $G-$pmorphism such that $\varsigma'\circ(\tau_W\circ \sigma_Y)=t$ and $\varsigma'\circ(\tau_W\circ \sigma_Z)=z$, then
\begin{center}
$\varsigma'(\overline{w})=\varsigma'(\tau_W(\sigma_Y(w)))=t(w)=\varsigma(\overline{w})$, if $w\in Y$,
\end{center}
\begin{center}
$\varsigma'(\overline{w})=\varsigma'(\tau_W(\sigma_Z(w)))=s(w)=\varsigma(\overline{w})$, if $w\in Y$,
\end{center}
Thus, $\varsigma$ is unique.

$$
\xymatrix@=1.5cm{
X \ar[r]^{f} \ar[d]_{g} & Y\ar[d]^{\tau_W\circ \sigma_Y} \ar@/_{-2pc}/[rdd]^{t} &\\
Z\ar[r]_{\tau_W\circ \sigma_Z}  \ar@/_1pc/[drr]_{s}   & W/\varepsilon_W\ar@{-->}[rd]^{\varsigma}&\\
  &   &R
  }
$$
\end{proof}

\section[Enveloping Action of a Topological Partial Action]{Enveloping Action of a Topological Partial  Action}

In this section, we prove that passing from a partial action to its enveloping action can be seen as a functor $E$ between the categories $G-\mathbf{pActTop}$ and $G-\mathbf{ActTop}$.
Moreover, we prove that $E$ is a left adjoint of the inclusion functor $I$, which implies that it preserves co-products.

\vspace{0.4cm}

Recall that any topological partial action can be seen as the restriction of a topological global action.
We start by presenting some important definitions and propositions, which will be used in the categorical context that we want.\\

\begin{definition}[\cite{Abadie2}]
Let $\alpha $ and $\beta$ be topological partial actions of $ G $ on $ X $ and $ Y,$ respectively. It is said that $ \alpha $ and $ \beta $ are $G-$pequivalent, if there exists a $G-$pequivalence $f:X\to Y$. That is, $f$ is a homeomorphism such that $f(X_g)\subseteq Y_g$ and $\beta_g(f(x))=f(\alpha_g(x))$, for each  pair $g\in G$ and $x\in X_{g^{-1}}$.\\
\end{definition}

\begin{definition}
Let $ \alpha $ be a topological partial action  of $ G $ on $ X $. We  say that $ (Y, \psi,j) $ is a  globalization of $ \alpha $, if $ \psi $ is a topological global  action of $ G $ on $ Y $ and   $ j: X \to Y $ is an injective continuous function such that $ j (X) $ is an open subset of $Y$ and the induced topological partial action of $ G $ on $ j (X) $ and $ \alpha $ are $G-$pequivalent.\\
\end{definition}

Given a topological partial action $\alpha$, the following propositions show the existence of a minimal and unique globalization for $\alpha$, which is called the enveloping action.\\

\begin{proposition}[\cite{Abadie2}]\label{equi}
Let $ \alpha $ be a topological partial action  of $ G $ on $ X $. The relation $ \sim $ on $ G \times X $ given by $ (g, x) \sim (h, y) $ if, and only if, $ x \in X_{g^{- 1}h} $ and $ \alpha_ {h^{- 1}g} (x) = y $, is an equivalence relation.\\
\end{proposition}

\begin{theorem}[\cite{Abadie2}]\label{indc}
Let $ \sim $ the equivalence relation given in the previous proposition and $ X^{\mathbf {e}} = (G \times X) / \sim $.
Then:
\begin{enumerate}
  \item The quotient function $ q_X: G \times X \to X^\mathbf{e} $, is open.
  \item The function $ \alpha^\mathbf{e}: G \times X^\mathbf{e} \to X^\mathbf{e} $ given by $ \alpha^\mathbf{e} (g, q_X (h, x)) = q_X (gh, x),$ is a topological global  action of $ G $ on $ X^\mathbf{e}$.
  \item The function $ \iota_X: X \to X^\mathbf{e} $ given by $ \iota_X(x) = q _X(e, x ) $, for all $x\in X$, is injective
  and open. In particular, $\iota_X(X) $ is an open subset of $ X^\mathbf{e} $.
  \item $ (X^\mathbf{e}, \alpha^\mathbf{e},\iota_X) $ is a globalization of $ \alpha $.
\end{enumerate}
\end{theorem}
The details of the proof can be found in \cite{Avila-Molina}. In the same work, the following notions and results are explicitly introduced.\\

\begin{definition}[\cite{Avila-Molina}]
Let $ \alpha $ be a topological partial action  of $ G $ on $ X $. We say that the globalization $ (Z, \phi, m) $ of $ \alpha $ is minimal, if for every globalization $ (Y, \psi, j) $ of $ \alpha $, there exists a unique $ G-$morphism $ \mu: Z \to Y $ such that $ j = \mu \circ m $.\\
\end{definition}
In the last definition, the existence of $ \mu $ guarantees that the following diagram commutes:
$$
\xymatrix@=1.5cm{
X \ar[dr]_j \ar[r]^m & Z\ar@{-->}[d]^{\mu}\\
& Y
}
$$

\begin{definition}[\cite{Avila-Molina}]
Let $ \varphi $ and $ \psi $ be topological global actions of $ G $ on $ X $ and $ Y$, respectively. We say that $ \varphi $ and $ \psi $ are $G-$pequivalent, if there exists a homeomorphism $ f: X \to Y $  such that $f(\varphi(g,x))=\psi(g,f(x))$ for each pair $g\in G$ and $x\in X$.\\
\end{definition}

\begin{proposition}[\cite{Avila-Molina}]\label{uni}
Let $ \varphi $ and $ \psi $ be topological global actions of $ G $ on $ X $. If $ (Z, \phi, m) $ is a minimal globalization of $ \alpha $, then $ (Z, \phi, m)$  is unique, except for $G-$pequivalences.\\
\end{proposition}

Although the following theorem was proposed in \cite{Abadie2}, the details of
its proof can be found in \cite{Avila-Molina}.\\

\begin{theorem}[\cite{Abadie2}]\label{Teo4}
If $ \alpha $ is a  topological partial action of $ G $ on $ X $, then with the notations of Theorem \ref{indc}, $ (X^\mathbf{e}, \alpha^\mathbf{e}, \iota_X) $ is a minimal globalization of $ \alpha $.\\
\end{theorem}

\begin{definition}[\cite{Abadie2}]
Let $ \alpha $ be a topological partial  action of $ G $ on $ X $. The minimal  globalization of $\alpha$,  given by $ (X^\mathbf{e}, \alpha^\mathbf{e}, \iota_X)$,  is called the enveloping  of $ \alpha $, $X^\mathbf{e}$ is  called the enveloping space of $X$, and $\alpha^\mathbf{e}$ is called the enveloping action of $\alpha$.\\
\end{definition}

In what follows, we will prove that passing to the enveloping action can be seen as a functor between the categories $G-\mathbf{pActTop}$ and $G-\mathbf{ActTop}$, which is a left adjoint of the inclusion functor.
For this,  it is necessary to define $E:(G-\mathbf{pActTop})\to (G-\mathbf{ActTop})$ as follows and prove that it is a covariant functor.
\begin{enumerate}
  \item If $X\in \text{Ob}(G-\mathbf{pActTop})$, then $E(X)=X^\mathbf{e}$, where $X^\mathbf{e}$ is the enveloping space of $X$.
Thus, $X^\mathbf{e}\in\text{Ob}(G-\mathbf{ActTop})$.
  \item If $X,Y\in \text{Ob}(G-\mathbf{pActTop})$ and  $f\in\text{Mor}(X,Y)$, then $E(f)=f^\mathbf{e}$, where $f^\mathbf{e}:X^\mathbf{e}\to Y^\mathbf{e}$ is given by $f^\mathbf{e}(q_X(h,x))=q_Y(h,f(x))$.
Then $f^e$ is a  $G-$morphism.
Indeed,
      \begin{enumerate}
        \item Suppose that $ \alpha $ and $ \beta $ are topological partial actions of $ G $ on $ X $ and $ Y,$ respectively.
If $ q_X (g, x) = q_X (h, y)$,  then $x\in X_{g^{-1}h}$ and $\alpha_{h^{-1}g}(x)=y$.
Since $f$ is a  $G-$pmorphism, then $f(x)\in Y_{g^{-1}h}$ and $\beta_{h^{-1}g}(f(x))=f(\alpha_{h^{-1}g}(x))=f(y).$ So
            \begin{align*}
            f^\mathbf{e}(q_X(g,x))&=q_Y(g,f(x))=q_Y(h,f(y))=f^\mathbf{e}(q_X(h,y)).
            \end{align*}
            Thus, $f^\mathbf{e}$ is well defined.
        \item Given  $q_X(h,x)\in X^\mathbf{e}$ and an open set $U\in Y^\mathbf{e}$ such that $f^\mathbf{e}(q_X(h,x))\in U$, we have that  $q_Y(h,f(x))\in U$.
Since $ q_Y$ is continuous, $ q ^ {- 1} _Y (U) $ is an open neighborhood of $ (h, f (x)) $ and therefore there are open neighborhoods $ M $ and $ N $ of $ h $ and $ f (x) $ respectively such that  $M\times N\subseteq q^{-1}_Y(U)$. Now, since $ f $ is continuous, $ f^{- 1}(N) $ is an open neighborhood of $ x $. Thus, $q_X(M\times f^{-1}(N))$ is an open neighborhood of $ q_X (h, x) $. Furthermore, $f^\mathbf{e}(q_X(M\times f^{-1}(N)))\subseteq U$.
Hence, $ f^\mathbf{e} $ is continuous.
 \item Now, if $q_X(h,x)\in X^\mathbf{e}$ and $g\in G$, then
\begin{align*}
f^\mathbf{e}(\alpha^\mathbf{e}(g,q_X(h,x)))&=f^\mathbf{e}(q_X(gh,x))\\
&=q_Y(gh,f(x))\\
&=\beta^\mathbf{e}(g,q_Y(h,f(x))\\
&=\beta^\mathbf{e}(g,f^\mathbf{e}(q_X(h,x))).
\end{align*}
\end{enumerate}
Consequently, $f^\mathbf{e}$  is an element of $\text{Mor}(X^\mathbf{e},Y^\mathbf{e})$.
\item  If $X,Y,Z\in \text{Ob}(G-\mathbf{pActTop})$, $f_1\in\text{Mor}(X,Y)$ and $f_2\in\text{Mor}(Y,Z)$, then for each $q_X(g,x)\in X^\mathbf{e}$ we have that
\begin{equation*}
\begin{split}
E(f_2\circ f_1)(q_X(g,x))&=(f_2\circ f_1)^\mathbf{e}(q_X(g,x))\\
&=q_Z(g,(f_2\circ f_1)(x))\\
&=q_Z(g,f_2(f_1(x)))\\
&=f_2^\mathbf{e}(q_Y(g,f_1(x)))\\
&=f_2^\mathbf{e}(f_1^\mathbf{e}(q_X(g,x))\\
&=E(f_2)(E(f_1)(q_X(g,x)))\\
&=(E(f_2)\circ E(f_1))(q_X(g,x)).
\end{split}
\end{equation*}
\item  For each $q_X(g,x)\in X^\mathbf{e}$ we have
\begin{align*}
E(i_{X})(q_X(g,x))
&=i_{X}^\mathbf{e}(q_X(g,x))\\
&=q_X(g,i_{X}(x))\\
&=q_X(g,x)\\
&=i_{X^\mathbf{e}}(q_X(g,x))\\
&=i_{E(X)}(q_X(g,x)).
\end{align*}
\end{enumerate}

Thus, we have obtained de following result.\\

\begin{theorem}
Let the function $E:(G-\mathbf{pActTop})\to (G-\mathbf{ActTop})$ defined as $E(X)=X^e$ for every $X\in Ob(G-\mathbf{pActTop})$ and $E(f)=f^e$ for every $f\in Mor(X,Y)$ with $X,Y\in Ob(G-\mathbf{pActTop})$ and where   $f^\mathbf{e}:X^\mathbf{e}\to Y^\mathbf{e}$ is given by $f^\mathbf{e}(q_X(h,x))=q_Y(h,f(x))$. Then $E$ defines a covariant functor between the categories $G-\mathbf{pActTop}$ and $G-\mathbf{ActTop}$. \\
\end{theorem}

The functor $E$ will be called the enveloping functor from $G-\mathbf{pActTop}$ to $G-\mathbf{ActTop}$. To conclude the study of $E$, it is necessary to recall the following notions.\\

\begin{definition}
Let $ \mathbf{C} $ and $ \mathbf{D} $ be categories and $ F, K $  covariant functors from $ \mathbf{C} $ to $ \mathbf {D} $. A function $ \theta: F \to K $ is said to be a natural transformation if it associates to every object $ X $ of $ \mathbf{C}$ a morphism $ \theta_X: F (X) \to K (X) $ in $ \mathbf{D}$, such that for every morphism $ f: X \to Y $ in $ \mathbf {C} $ it holds $\theta_{Y}\circ F(f) = K(f)\circ\theta_X$. That is, the following diagram is commutative:
$$
\xymatrix@=1.5cm{
F(X) \ar[r]^{\theta_X} \ar[d]_{F(f)} & K(X)\ar[d]^{K(f)}\\
F(Y)\ar[r]_{\theta_Y}& K(Y)
}
$$
\end{definition}

\begin{proposition}\label{nat}
If $ \mathbf{C} = G-\mathbf{pActTop} $, $ \mathbf{D} = G-\mathbf{ActTop} $, $ E $ is the enveloping functor from $ \mathbf{C} $ to $ \mathbf {D} $, and $ I $ is the inclusion functor from $ \mathbf{D} $ to $ \mathbf{C} $, then $ \theta: i_{\mathbf {C}} \to IE $ where $ \theta_X = \iota_X: X \to X^\mathbf{e} $, is a natural transformation.
\end{proposition}
\begin{proof}
If $ X \in \text{Ob} (\mathbf{C}) $, then there exists a topological partial action $ \alpha $ of $ G $ on $ X $ such that $ \theta_X = \iota_X: X \to X^\mathbf{e} $ is a $G-$pmorphism.
Now, if $ f \in \text{Mor} (X, Y) $, then there exists a partial topological action $ \beta $ of $ G $ on $ Y $, and therefore $ \theta_Y = \iota_Y: Y \to Y^\mathbf {e} $ is a $G-$pmorphism such that for each $ x \in X $ we have
\begin{align*}
(IE(f)\circ \theta_X)(x)&= (IE(f))(\theta_X(x))\\
&= I(E(f))(\iota_X(x))\\
&= E(f)(q_X(e,x))\\
&=q_Y(e,f(x))\\
&=\iota_Y(f(x))\\
&=\theta_Y(f(x))\\
&=\theta_Y(i_{\mathbf{C}}(f)(x))\\
&=(\theta_Y\circ i_{\mathbf{C}}(f))(x).
\end{align*}
Thus,  $IE(f)\circ \theta_X= \theta_Y\circ i_{\mathbf{C}}(f)$.
$$
\xymatrix@=1.5cm{
i_{\mathbf{C}}(X) \ar[r]^{\theta_X} \ar[d]_{i_{\mathbf{C}}(f)} & IE(X)\ar[d]^{IE(f)}\\
i_{\mathbf{C}}(Y)\ar[r]_{\theta_Y}& IE(Y)
}
$$
\end{proof}

\begin{lemma}[\cite{Stonek}]
Let $ \mathbf{C} $ and $ \mathbf{D} $ be  categories, $ F $  a covariant functor from $ \mathbf {C} $ to $ \mathbf {D} $ and $ K $ be a covariant functor from $ \mathbf{D} $ to $ \mathbf{C} $. The pair $ (F, K) $ is adjoint if, and only if, there exists a natural transformation $ \theta: i _ {\mathbf {C}} \to KF $ such that for every $ X \in \text {Ob}(\mathbf{C}) $, $ Y \in \text{Ob} (\mathbf{D}) $ and $ f \in \text{Mor} (X, K (Y)) $, there exists a unique $ g \in \text{Mor} (F (X), Y) $ such that $ f = K (g) \circ \theta_X $. That is, the following diagram is commutative:
$$
\xymatrix@=1.5cm{
X \ar[r]^f \ar[d]_{\theta_X} & K(Y)\\
KF(X)\ar@{-->}[ur]_{K(g)}&
}
$$
\end{lemma}

\begin{proposition}[\cite{Abadie2}]\label{propuni}
Let  $ \alpha $ be  a topological partial action of $ G $ on $ X $. If $\beta $ is a topological global action of $ G $ on $ Y $, then $ \beta $ is a partial topological action of $ G $ on $ Y $ such that for every $G-$pmorphism $ f: X \to Y $, there exists a unique $G-$morphism $ \tilde{f}: X^\mathbf{e} \to Y $ such that $ f = \tilde{f} \circ \iota_X $. That is, the following diagram is commutative:
$$
\xymatrix@=1.5cm{
X \ar[r]^f \ar[d]_{\iota_X} & Y\\
X^\mathbf{e}\ar@{-->}[ur]_{\tilde{f}}&
}
$$
\end{proposition}
\begin{proof}
It is clear that $\beta$ is a partial topological action of $ G $ on $ Y $.
If $ f: X \to Y $ is a $G-$pmorphism, the function $ \tilde{f}: X^\mathbf {e} \to Y $ given by $$ \tilde {f} (q_X (g, x)) = \beta(g, f (x))$$ is well defined.
Indeed, if $ q_X (g, x) = q_X (h, y) $, we have that $ x \in X_{g^{- 1}h} $ and $ \alpha_ {h^{- 1}g} (x) = y $.
Thus, $ \beta_{h^{- 1} g} (f (x)) = f(\alpha_{h^{- 1} g} (x)) = f (y) $ and consequently
\begin{align*}
(\beta_{h^{-1}}\circ\beta_{g})(f(x))&=\beta(h^{-1},\beta(g,f(x)))\\
&=\beta(h^{-1}g,f(x))\\
&=\beta_{h^{-1}g}(f(x))\\
&=f(y).
\end{align*}
Hence, $\beta_{g}(f(x))=\beta_{h}(f(y))$ and so
\begin{align*}
\tilde{f}(q_X(g,x))=\beta(g,f(x))=\beta(h,f(y))=\tilde{f}(q_X(h,y)).
\end{align*}

Now, suppose that $q_X(g,x)\in X^\mathbf{e}$ and $U\subseteq Y$ is an open neighborhood of $\tilde{f}(q_X(g,x))=\beta(g,f(x))$.
Since $ \beta $ is continuous, there are open sets $ M \subseteq G $ and $ N \subseteq Y $ such that $(g,f(x))\in M\times N$ and $\beta(M\times N)\subseteq U$.
Since the function $ f $ is continuous, there is an open set $ S \subseteq X $ such that $ x \in S $ and $ f (S) \subseteq N $.
Also, since the quotient function $ q_X $ is open, $V=q_X(M\times S)\subseteq X^\mathbf{e}$ is an open neighborhood of $ q_X (g, x) $ such that  $\tilde{f}(V)=\tilde{f}(q_X(M\times S))=\beta(M\times f(S))\subseteq\beta(M\times N)\subseteq U$.
Therefore, $ \tilde{f} $ is continuous.
Also, if $ g, h \in G $ and $ x \in X $, then
\begin{align*}
\tilde{f}(\alpha^\mathbf{e}(g,q_X(h,x)))
&=\tilde{f}(q_X(gh,x))\\
&=\beta(gh,f(x))\\
&=\beta(g,\beta(h,f(x)))\\
&=\beta(g,\tilde{f}(q_X(h,x))).
\end{align*}
Therefore, $ \tilde {f} $ is a $ G-$morphism. Finally,
\begin{align*}
(\tilde{f}\circ \iota_X)(x)=\tilde{f}(\iota_X(x))=\tilde{f}(q_X(e,x))=\beta_e(f(x))=f(x),
\end{align*}
for each $x\in X$.
Furthermore, if $\hat{f}:X^\mathbf{e}\to Y$ is a $ G-$morphism such that $f=\hat{f}\circ \iota_X$, then for each $q_X(g,x)\in X^\mathbf{e}$ we have
\begin{align*}
\tilde{f}(q_X(g,x))
&=\beta(g,f(x))\\
&=\beta(g,\hat{f}(\iota_X(x)))\\
& =\beta(g,\hat{f}(q_X(e,x)))\\
&=\hat{f}(\alpha^\mathbf{e}(g,q_X(e,x)))\\
&=\hat{f}(q_X(g,x)).
\end{align*}
So $\tilde{f}$ is unique.
\end{proof}
\vspace{0.2cm}

\begin{theorem}\label{adjizq}
Let $E$ the enveloping functor from $\mathbf{C}=G-\mathbf{pActTop}$ to $\mathbf{D}=G-\mathbf{ActTop}$  and $I$
the inclusion functor from $\mathbf{D}$ to $\mathbf{C}$. Then $(E,I)$ is an adjoint pair and therefore, $ E $ is a left adjoint  of $ I $.
\end{theorem}
\begin{proof}
Let $X\in \text{Ob}(\mathbf{C})$, $Y\in \text{Ob}(\mathbf{D})$ and $f\in\text{Mor}(X,I(Y))$.  Then there exists a partial topological action $ \alpha $ of $ G $ on $ X $ and a topological global  action $\beta $ of $ G $ on $ Y $. Since $ I (Y) = Y $, then $\beta$ is a topological partial action of $ G $ on $ Y $ and so $ f: X \to Y $ is a $G-$pmorphism.\\

By Proposition \ref{nat}, $\theta: i_{\mathbf{C}} \to IE$, where $\theta_X=\iota_X:X\to X^\mathbf{e}$,
 is a natural transformation. Also, by Proposition \ref{propuni}, there exists a unique $ G-$morphism $\tilde{f}:X^\mathbf{e}\to Y$ such that $f=\tilde{f}\circ \iota_X$. So $\tilde{f}\in \text{Mor}(E(X),Y)$ is unique and $f=I(\tilde{f})\circ \theta_X$. That is, the following diagram is commutative.
$$
\xymatrix@=1.5cm{
X \ar[r]^f \ar[d]_{\theta_X} & I(Y)\\
IE(X)\ar@{-->}[ur]_{I(\tilde{f})}&
}
$$
Thus $(E,I)$ is an adjoint pair and therefore the enveloping functor $E$ is a left adjoint of the inclusion functor $I$.\\
\end{proof}
In category theory, it is well-known that if $\mathbf{C}$ and $\mathbf{D}$ are two categories, $ F $ is a covariant functor from $\mathbf{C}$  to \textbf{D} and $K$ is a covariant functor from $\mathbf{D}$ to $\mathbf{C}$ such that $ F $ is a left adjoint of $ K$, then $ F $ preserves co-products and $ K $ preserves products. Thus,  we have the following results.\\
\begin{corollary}
Let $ G $ be a topological group and $\{X_i\}_{i\in I}$ a collection of topological spaces. If $ \alpha_i $ is a topological partial action of $ G $ on $ X_i $ and $(X_i^\mathbf{e},\alpha^\mathbf{e}_i,\iota_i)$ is the enveloping of $ \alpha_i $ for every $ i \in I $, then the enveloping of the topological partial action of $ G $ on the co-product  $X=\bigsqcup_{i\in I}X_i$ is $(X^\mathbf{e},\alpha^\mathbf{e},\iota)$ where $X^\mathbf{e}=\bigsqcup_{i\in I}X_i^\mathbf{e}$, $\alpha^\mathbf{e}:G\times X^\mathbf{e}\to X^\mathbf{e}$ is the topological global  action given by $\alpha^\mathbf{e}(g,q_X(h,x))=\alpha^\mathbf{e}_{i}(g,q_X(h,x))$, for each pair $g\in G$ and $q_X(h,x)\in X_i^\mathbf{e}$, and $\iota:X\to X^\mathbf{e}$
is the function given by  $\iota(x)=\iota_{i}(x)$ for each $x\in X_i$.\\
\end{corollary}

\begin{corollary}
Let $ G $ be a topological group and $\{X_i\}_{i\in I}$ a collection of topological spaces. If $ \varphi_i $ is a topological global action of $ G $ on $ X_i $ for every $i\in I$ and $X=\prod_{i\in I}X_i$ is endowed with the product topology, then the product object of the collection $\{X_i\}_{i\in I}$ in the category $G-\mathbf{pActTop}$ is $X$ together with the topological global action $\varphi:G\times X \to X$  given by $\varphi(g,x)=(\varphi_i(g,x_i))_{i\in I}$, where $g\in G$ and $x\in X$. That is, the object product of the collection $\{X_i\}_{i\in I}$ in the category $G-\mathbf{ActTop}$, given in Theorem \ref{proglo}, coincides with object product of the collection $\{X_i\}_{i\in I}$ in the category $G-\mathbf{pActTop}$. 
\end{corollary}




\section*{Financial Support}
The first author received funding from the Universidad del Tolima. The second author received funding from the Universidad de Ibagué and is currently receiving funding from the Universidad de Oviedo.

\section*{Acknowledgments}


\end{document}